\documentclass[11pt, oneside]{amsart}

\usepackage{amsmath,amsthm,amssymb}  

\usepackage{graphicx}    
\graphicspath{ {images/} }  
\numberwithin{figure}{section}  

\usepackage{geometry}  

\usepackage{gensymb}
\usepackage{latexsym}
\usepackage{graphicx}
\usepackage{epsfig,epsf,rotating}
\usepackage{subfigure}
\usepackage{epstopdf}

\usepackage{tikz-cd}
\usepackage{amsthm} 
\usepackage[capitalise]{cleveref} 
\theoremstyle{plain}
\newtheorem{theorem}{Theorem}[section]
\newtheorem{corollary}[theorem]{Corollary}
\newtheorem{lemma}[theorem]{Lemma}
\newtheorem{proposition}[theorem]{Proposition}

\theoremstyle{definition}
\newtheorem{definition}[theorem]{Definition}
\theoremstyle{remark}

\newtheorem{remark}[theorem]{Remark}

\title{Tau invariants for balanced spatial graphs}
\author{Katherine Vance$^{\dag}$}
\address{Department of Mathematics \\Simpson College \\701 North C Street, Indianola, IA 50125}  
\date{\today}                                          

\thanks{$^{\dag}$ This work was partially supported National Science Foundation grant DMS-1309070.}

\begin{document}

\begin{abstract}
In 2003, Ozsv\'ath and Szab\'o defined the concordance invariant $\tau$ for knots in oriented 3-manifolds as part of the Heegaard Floer homology package.   In 2011, Sarkar gave a combinatorial definition of $\tau$ for knots in $S^3$ and a combinatorial proof that $\tau$ gives a lower bound for the slice genus of a knot.  Recently, Harvey and O'Donnol defined a relatively bigraded combinatorial Heegaard Floer homology theory for transverse spatial graphs in $S^3$, extending HFK for knots. We define a $\mathbb{Z}$-filtered chain complex for balanced spatial graphs whose associated graded chain complex has homology determined by Harvey and O'Donnol's graph Floer homology.   We use this to show that there is a well-defined $\tau$ invariant for balanced spatial graphs generalizing the $\tau$ knot concordance invariant.   In particular, this defines a $\tau$ invariant for links in $S^3$.   Using techniques similar to those  of Sarkar, we show that our $\tau$ invariant is an obstruction to a link being slice.
\end{abstract}

\maketitle

\section{Introduction}

\label{ch:Intro}

\subsection{Background}
 
	A graph is a one-dimensional CW-complex whose edges (one-cells) may be oriented.  A spatial graph is a smooth or piecewise linear embedding $f:G \rightarrow S^3$, where $G$ is an (oriented) graph.  One way to think of spatial graphs is as a generalization of the classical study of knots and links, which are embeddings of one or more ordered $S^1$ components into $S^3$.  Just as for knots and links, we consider spatial graphs up to ambient isotopy.  
	
\begin{figure}
 	\begin{center}
    	\includegraphics[width=.25\textwidth]{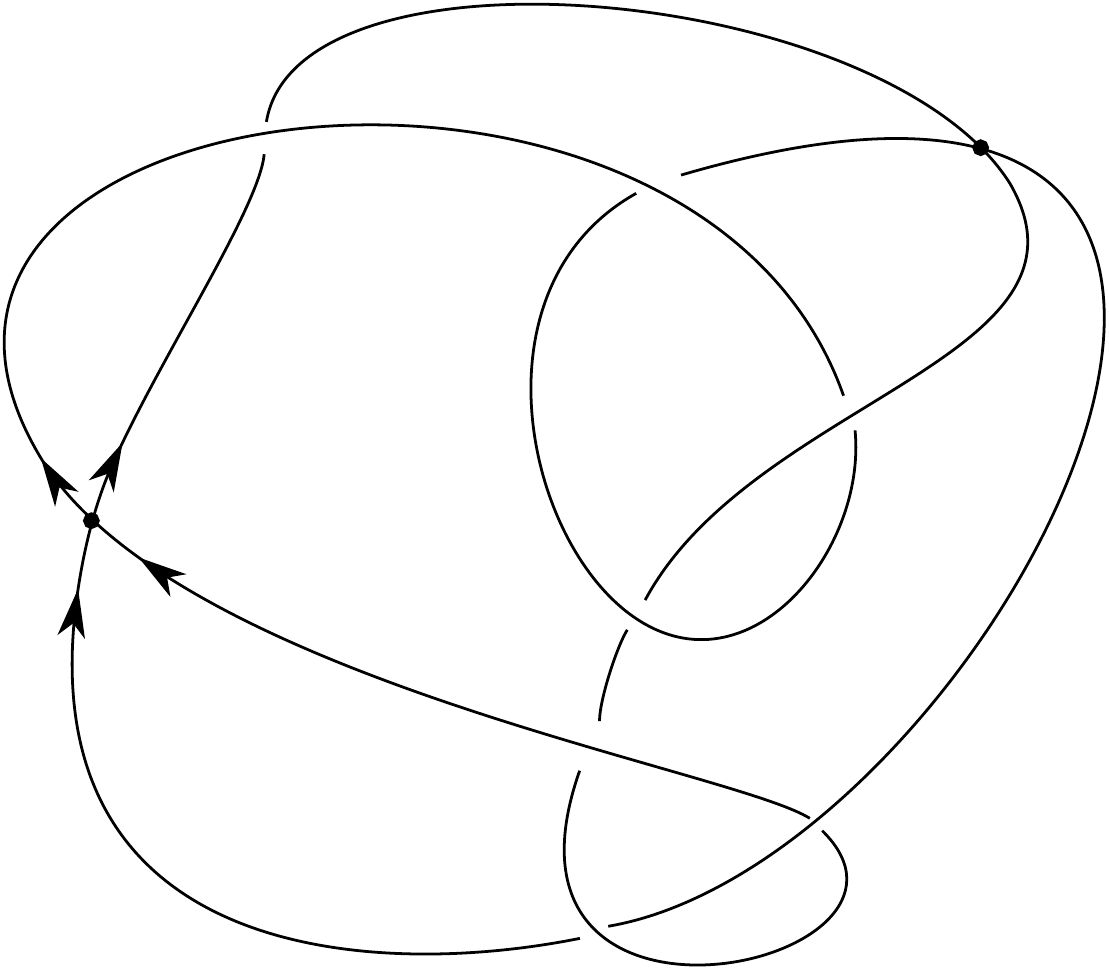}
  	\end{center}
  	\caption{A spatial graph}
\end{figure}
	
	Knot Floer homology is a package of invariants which was independently defined in 2002 by Ozsv\'ath and Szab\'o \cite{MR2065507} and by Rasmussen \cite{Rasmussen}.  One invariant from the knot Floer homology package is the $\tau$ invariant, which was defined by Ozsv\'ath and Szab\'o in 2004 \cite{MR2023281}.   

	One reason the $\tau$ invariant is important is its relationship to knot concordance.  The $\tau$ invariant is a concordance invariant and its absolute value is a lower bound for slice genus \cite{MR2065507}.  In 2011, Sarkar gave a combinatorial proof of the relationship between $\tau$ and slice genus \cite{MR2915478}.   Recently, Harvey and O'Donnol have defined graph Floer homology for a certain class of spatial graphs in $S^3$ using a grid diagram construction analogous to that used for knots and links \cite{Harvey-ODonnol}.  However, while knot Floer homology is filtered by the integers, Harvey and O'Donnol's graph Floer homology is not; rather it is relatively graded graded by the first homology group of the spatial graph complement.  

\subsection{Summary of main results}
In this paper, we define a filtered version of graph Floer homology for balanced transverse spatial graphs whose associated graded object is Harvey and O'Donnol's $HFG$ and prove that it is a spatial graph invariant.  We prove that the filtered graph Floer chain complex is, up to filtered quasi-isomorphism, an invariant of balanced spatial graphs.  Thus we have the following theorem.

{
\renewcommand{\thetheorem}{\ref{thm:main}}
\begin{theorem}
For grid diagrams $g, g'$ representing $f:G \rightarrow S^3$,  there exist filtered quasi-isomorphisms  $\phi_1: CF^-(g) \rightarrow CF^-(g')$ and $\phi_2: CF^-(g') \rightarrow CF^-(g)$ which preserve the symmetrized filtration $\lbrace \mathcal{F}^{-H}_s \rbrace$.
\end{theorem}
\addtocounter{theorem}{-1}
}

This allows us to define a $\tau$ invariant for balanced spatial graphs and prove that it is an invariant.

{
\renewcommand{\thetheorem}{\ref{def:tau}}
\begin{definition} For a graph grid diagram $g$ representing a balanced spatial graph $f:G \rightarrow S^3$, define 	the $\tau$ invariant of $g$ to be	\[ \tau (g) = \text{min} \{ m \in \frac 1 2 \mathbb{Z} | \iota _m \text{ is non-trivial} \}\]
	where $\iota_m : H_*(\widehat{\mathcal{F}}_m^H) \rightarrow H_*(\widehat{CF}(g))$ is the map induced by inclusion.
\end{definition}
\addtocounter{theorem}{-1}
}

{
\renewcommand{\thetheorem}{\ref{cor:tau}}
\begin{corollary}
	If $g$ and $\overline{g}$ are graph grid diagrams representing a balanced spatial graph $f: G \rightarrow S^3$, then $\tau (g) = \tau (\overline{g})$.
\end{corollary}
\addtocounter{theorem}{-1}
}

Considering links as spatial graphs with one vertex and one edge in each link component, we obtain the following result relating the $\tau$ invariant to link cobordisms.

{
\renewcommand{\thetheorem}{\ref{thm:LinkCob}}\begin{theorem}
If $L_1$ and $L_2$ are $l_1$- and $l_2$-component links, respectively, and $F$ is a connected genus $g$ cobordism from $L_1$ to $L_2$, then 
\[ 1 - g - l_1 \leq \tau (L_1) - \tau (L_2) \leq g + l_2 -1. \]
\end{theorem}
\addtocounter{theorem}{-1}
}

As a corollary, we see that the $\tau$ invariant can be an obstruction to a link being slice.

{
\renewcommand{\thetheorem}{\ref{cor:linkslice}}
\begin{corollary}
If an $l$-component link $L$ has $\tau (L) > 0$ or $\tau (L) \leq -l$, then $L$ is not slice.
\end{corollary}
\addtocounter{theorem}{-1}
}

Recently, Cavallo independently defined a $\tau$ invariant for links and proved a result similar to \cref{thm:LinkCob} \cite{Cavallo}.

\section{Graph Floer Homology}
\label{ch:HFG}

In this section we give an overview of Harvey and O'Donnol's graph Floer homology, which is defined for transverse spatial graphs.  For precise definitions of spatial graphs and transverse spatial graphs, see \cite{Harvey-ODonnol}. 

\begin{figure}
	\begin{center}
	\includegraphics[width=.25\textwidth]{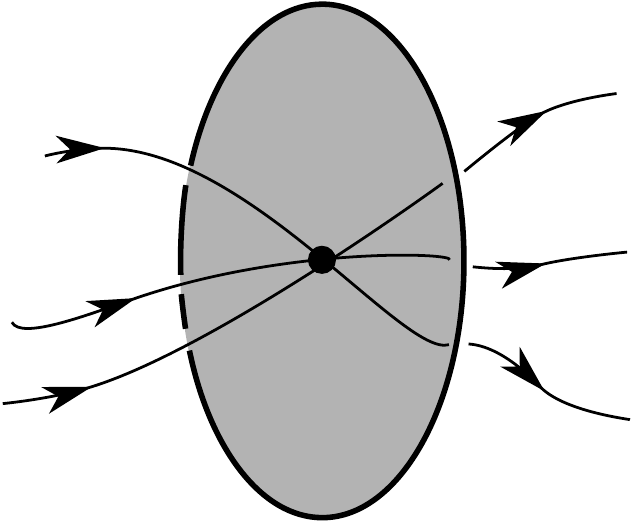}
	\caption{The standard disk separating incoming and outgoing edges at a vertex of a transverse spatial graph}
	\label{fig:transversedisk}
	\end{center}
\end{figure}

\begin{definition}
A \emph{spatial} graph is an embedding $f: G \rightarrow S^3$ of a $1$-dimensional CW-complex $G$ into $S^3$.  An \emph{oriented spatial graph} is a spatial graph with an orientation given for each edge.  For each vertex $v$ of an oriented spatial graph, the \emph{incoming edges} of $v$ are the edges incident to $v$ whose orientation points toward $v$, and the \emph{outgoing edges} of $v$ are the edges incident to $v$ whose orientation points away from $v$.  A \emph{disk graph} is one which has a standard disk $\mathcal{D}$ at each vertex, attached to the graph by identifying the center point of $\mathcal{D}$ with the vertex.  A \emph{transverse spatial graph} is an embedding $f: G \rightarrow S^3$ of an oriented disk graph $G$, such that at each vertex the standard disk is embedded in a plane that separates the incoming and outgoing edges, as shown in \cref{fig:transversedisk}.
\end{definition}

In contrast to spatial graph ambient isotopy, in which any combination of edges incident to a vertex can move freely, ambient isotopy of transverse spatial graphs only allows free movement of incoming edges with other incoming edges or outgoing edges with other outgoing edges at each vertex.  This is because the edges may not pass through the standard disk at the vertex.

Graph Floer homology is defined using grid diagrams, like the combinatorial definition of knot Floer homology.  The definition of spatial graph grid diagrams is very similar to the definition of grid diagrams for knots and links.  

	An index $n$ \textbf{graph grid diagram} for a transverse spatial graph is an $n$ by $n$ grid in which each grid square may contain an $O$-marking, an $X$-marking, or be empty, such that there is exactly one $O$ in each row and in each column.   We make a distinction between standard $O$-markings, which are those which are in the interior of a graph edge when we recover the spatial graph from the graph grid diagram, and special $O$-markings, which are vertices of the graph when it is recovered from the graph grid diagram.  We mark special $O$'s with an asterisk in the graph grid diagram.  Standard $O$-markings have exactly one $X$ in their row and column, while vertex $O$'s may have any number of $X$-markings in their row and column.  If a transverse spatial graph has more than one connected component, we require that there be at least one special $O$-marking in each component.  A \textbf{toroidal graph grid diagram} is one in which we think of the grid as being a torus, with the leftmost and rightmost gridlines identified and the top and bottom gridlines identified.

\begin{figure}
	\begin{center}
		\subfigure{\includegraphics[width=.25\textwidth]{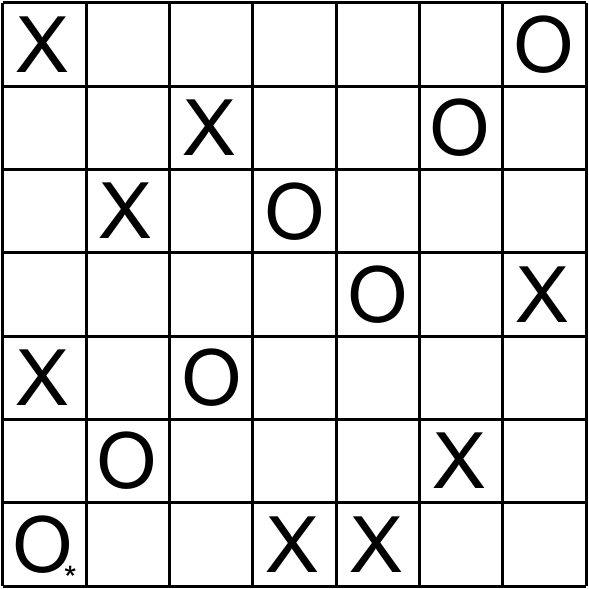}}
		\subfigure{\includegraphics[width=.25\textwidth]{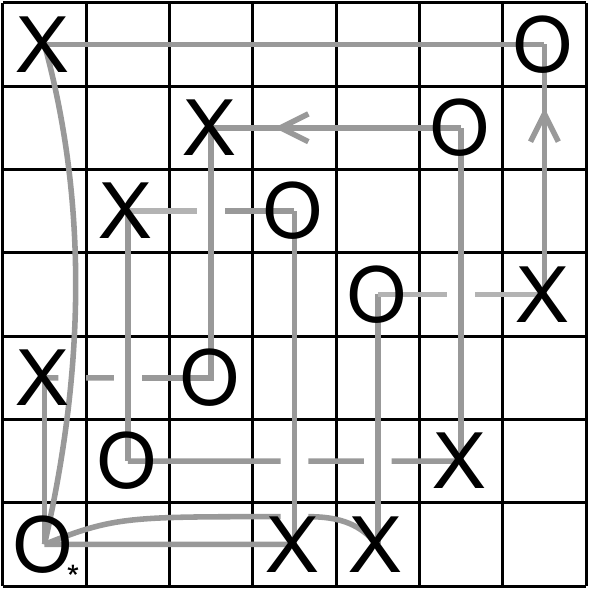}}
	\end{center}
	\caption{Graph grid diagram (note the starred vertex $O$)}
	\label{fig:gridknottedwedge}
\end{figure}

To recover the spatial graph from a grid diagram, connect the $X$'s to the $O$'s vertically and the $O$'s to the $X$'s horizontally.  At each crossing, the vertical strand is the overpass and the horizontal strand is the underpass. At vertex $O$'s (those with more than one $X$ in their row or column) use a straight line to connect the closest $X$ in the row or column to the vertex $O$ and a curved line to connect the more distant $X$'s to the vertex $O$, observing the same conventions with regard to the crossings created, so that the line connecting two markings within a column is always the overstrand.  See \cref{fig:gridknottedwedge}.  Just as is the case for knots and links, every transverse spatial graph can be represented by a graph grid diagram.

For knots and links, Cromwell's theorem \cite{MR1339757} gives a sequence of grid moves connecting any two grid diagrams representing equivalent links.  Harvey and O'Donnol have proved a similar theorem for transverse spatial graphs.

\begin{theorem}[\cite{Harvey-ODonnol}] 
\label{thm:gridmoves}
	Any two graph grid diagrams for a given transverse spatial graph are related by a finite sequence of cyclic permutation, commutation', and (de-)stabilization' moves.
\end{theorem}

\begin{figure}
	\begin{center}
		\subfigure{\includegraphics[width=.25\textwidth]{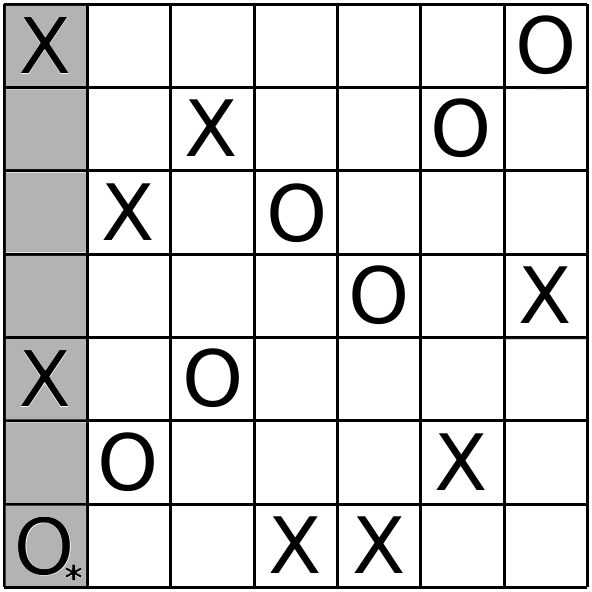}}
		\subfigure{\includegraphics[width=.25\textwidth]{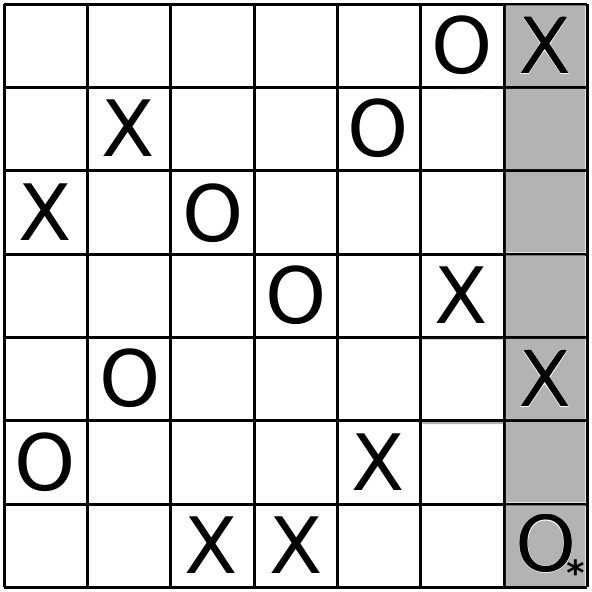}}
	\end{center}
	\caption{A cyclic permutation move}
	\label{fig:cyclicpermutation}
\end{figure}

A \textbf{cyclic permutation} moves the top (resp. bottom) row of a grid diagram to the bottom (resp. top) or moves the left (right) column  to the far right (left) of the diagram. See  the example in \cref{fig:cyclicpermutation}.  Thinking of the grid as a torus, this equates to changing which gridline we ``cut" the torus along to get the square diagram.

\begin{figure}
	\begin{center}
		\subfigure{\includegraphics[width=.25\textwidth]{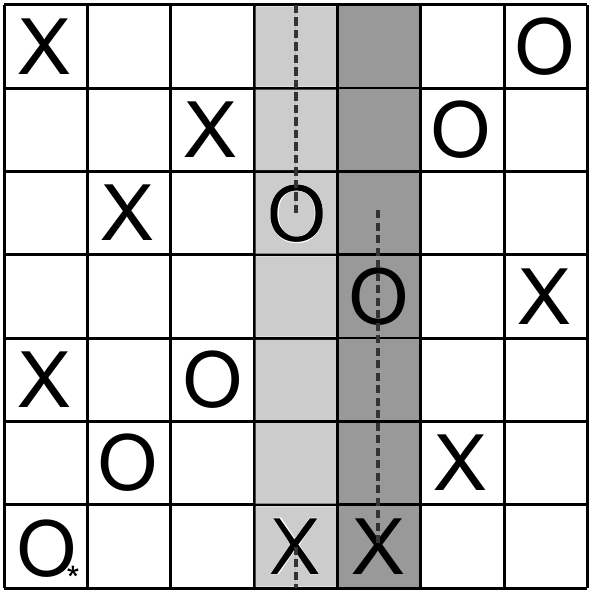}}
		\subfigure{\includegraphics[width=.25\textwidth]{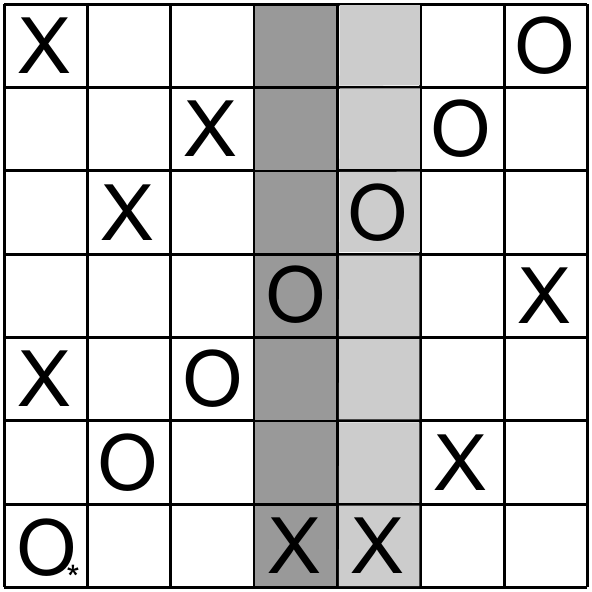}}
	\end{center}
	\caption{A commutation' move: notice the dotted helper arcs $LS_1, \, LS_2$ in the left-hand grid}
	\label{fig:gridcomm}
\end{figure}

Two adjacent columns (or rows) may be exchanged using a \textbf{commutation'} move if there are vertical (horizontal) line segments $LS_1$ and $LS_2$ on the torus such that $LS_1 \cup LS_2$ contain all the $X$'s and $O$'s in the two adjacent columns (rows), the projection of $LS_1 \cup LS_2$ to a single vertical circle $\beta_i$ (horizontal circle $\alpha_i$) is $\beta_i$ ($\alpha_i$), and the projection of their endpoints, $\partial(LS_1) \cup \partial(LS_2)$, to a single $\beta_i$ ($\alpha_i$) is precisely two points.    See the example in \cref{fig:gridcomm}.

\begin{figure}
	\begin{center}
		\subfigure{\includegraphics[width=.25\textwidth]{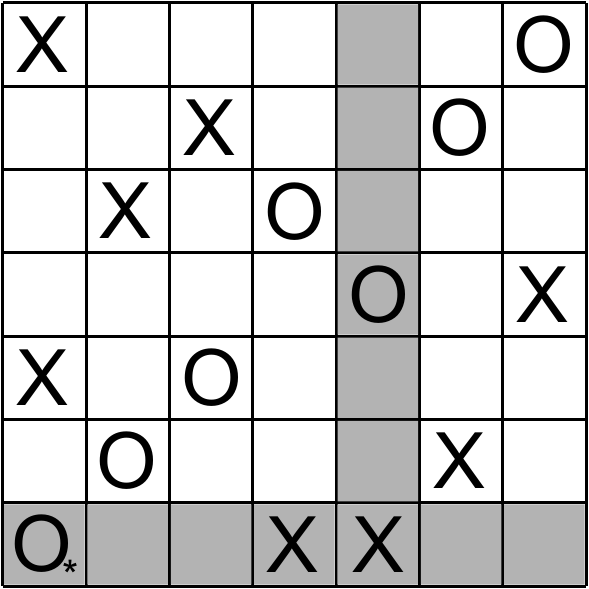}}
		\subfigure{\includegraphics[width=.29\textwidth]{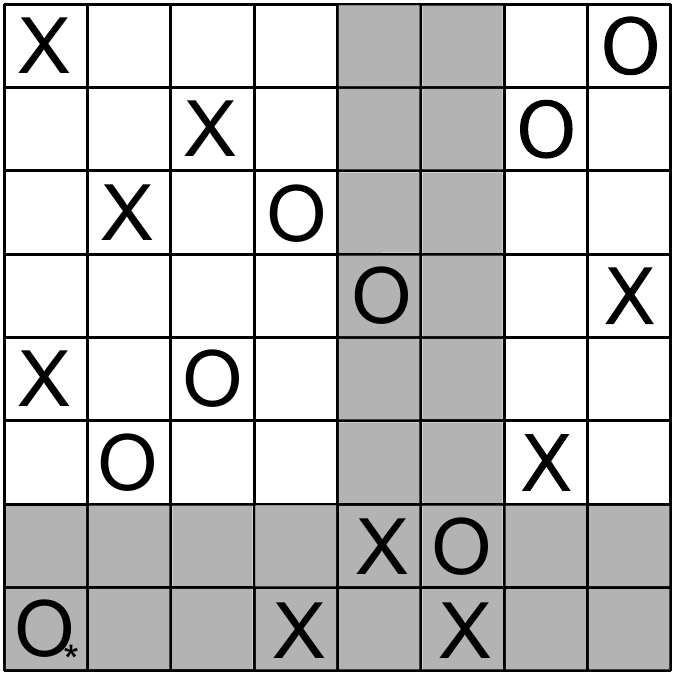}}
	\end{center}
	\caption{A row stabilization'}
	\label{fig:gridstab}
\end{figure}

A row (column) \textbf{stabilization'} at an $X$-marking is performed by adding one new row and one new column to the grid next to that $X$.  The $X$ is then moved to the new row (column), remaining in the same column (row), with the $O$ and any other $X$-markings in which were in the same row (column) as the $X$ being stabilized remaining in the old row (column). A new $X$-marking is placed in the intersection of the new column (row) and the row (column) previously occupied by the $X$-marking, and a new $O$ is placed in the intersection of the new row and column.  See the example in \cref{fig:gridstab}. A destabilization' is the opposite of a stabilization'.

Harvey and O'Donnol's graph Floer homology is defined for transverse spatial graphs without sinks or sources.  A \emph{sink} is a vertex with no outgoing edges and a \emph{source} is a vertex with no incoming edges.  In other words, graph Floer homology is defined for spatial graphs whose underlying graph has at least one incoming edge and at least one outgoing edge at every vertex.  This corresponds to a requirement that a graph grid diagram representing the spatial graph has at least one $X$-marking in every row and column.  

For a spatial graph $f: G \rightarrow S^3$ represented by an $n \times n$ graph grid diagram $g$, the graph Floer chain complex $(C^-(g), \partial^-)$ is freely generated as a module over $\mathbb{F}[U_1, ... , U_n]$, where $\mathbb{F} = \mathbb{Z}/2\mathbb{Z}$ and the $U_i$'s are formal variables corresponding to the $O$-markings $O_1, ... , O_n$ in the graph grid diagram.  The generating set of $C^-(g)$ is 
\[ \mathbf{S} = \lbrace \mathbf{x}=(x_1, ... , x_n) | x_i = \alpha_i \cap \beta_{\sigma(i)} \text{ for some } \sigma \in S_n \rbrace \]
where $S_n$ is the symmetric group on $n$ letters.

The map $\partial^- : C^-(g) \rightarrow C^-(g)$ counts empty rectangles in the toroidal graph grid diagram $g$.  An embedded rectangle $r$ in $g$ connects a generator $\mathbf{x}$ to another generator $\mathbf{y}$ if $x_i = y_i$ for all but two $i$, if $j < k$ are the two indices for which $\mathbf{x}$ and $\mathbf{y}$ are not equal, and if the corners of $r$ are, clockwise from the bottom left, $x_j, y_k, x_k,$ and $y_j$.  We say that $r$ is empty if the interior of $r$ does not contain any points of $\mathbf{x}$ or $\mathbf{y}$.  The set of empty rectangles from $\mathbf{x}$ to $\mathbf{y}$ is denoted $\mathcal{R}^\circ(\mathbf{x},\mathbf{y})$.  The map $\partial^- : C^-(g) \rightarrow C^-(g)$ is defined as follows on the generating set $\mathbf{S}$ and then extended to all of $C^-(g)$ as an $\mathbb{F}_2[U_1, ... , U_n]$-module homomorphism:
\[\partial^-(\mathbf{x}) = \sum_{\mathbf{y} \in \mathbf{S}} \sum_{\substack{r \in \mathcal{R}\degree(\mathbf{x},\mathbf{y}) \\ int(r) \cap \mathbb{X}= \varnothing}} U_1^{O_1(r)} \cdots U_n^{O_n(r)}\mathbf{y}\]
where $O_i(r)$ is zero if $O_i$ is not in $r$ and one if $O_i$ is in $r$.  Note that $\partial^-$ counts rectangles that contain any of the $O$-markings in $g$ but does not count any rectangles that contain $X$-markings.  This is because $H_1(S^3 - f(G))$ does not have a natural filtration, so Harvey and O'Donnol's graph Floer homology is graded rather than filtered.  

\begin{proposition}[\cite{Harvey-ODonnol} Proposition 4.10]
For $\partial^-:C^-(g) \rightarrow C^-(g)$ as defined above, $\partial^- \circ \partial^- = 0$.
\end{proposition}

\begin{figure}
 	\begin{center}
 	\def\svgwidth{0.25\textwidth}
    	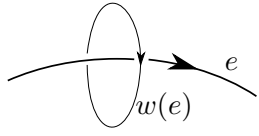
  	\end{center}
  	\caption{The weight of an edge}
  	\label{fig:edgeweight}
\end{figure}

Before we can define the Alexander grading we need to define weights of the edges of $G$.  We define a weight function $w: E(G) \rightarrow H_1(S^3 - f(G))$, where $E(G)$ is the set of edges of $G$, by mapping each edge $e \in E(G)$ to the homology class of the meridian of $e$, oriented according to the right-hand rule, as shown in \cref{fig:edgeweight}.

For $X$-markings and $O$-markings associated to the interior of an edge $e$, the weights are $w(X)=w(e)$ or $w(O) = w(e)$.  For $O$-markings associated to a graph vertex $v$, the weight is $w(O) = \sum_{e \in In(v)} w(e) = 
\sum_{e \in Out(v)} w(e)$, where $In(v)$ and $Out(v)$ are, respectively, the sets of incoming and outgoing edges of $v$.  

We can now define the Alexander grading on the generating set $\mathbf{S}$:
\[A(\mathbf{x}) = \sum_{p \in \mathbb{X}} \mathcal{J}(\mathbf{x},p)w(p) - \sum_{p \in \mathbb{O}} \mathcal{J}(\mathbf{x},p)w(p).\]
This grading is not well-defined on toroidal graph grid diagrams, but Harvey and O'Donnol show that the relative grading $A^{rel}(\mathbf{x}, \mathbf{y})= A(\mathbf{x}) - A(\mathbf{y})$ is well-defined on toroidal graph grid diagrams (\cite{Harvey-ODonnol} Corollary 4.14).

The graph Floer chain complex $(C^-(g), \partial^-)$ is bigraded, with an absolute $\mathbb{Z}$-valued grading (the Maslov grading) and a relative $H_1(S^3 - f(G))$-valued grading (the Alexander grading).  The graph Floer homology is $HFG^-(f) = H_*(C^-(g),\partial^-)$ for any graph grid diagram $g$ representing $f$, and it is also absolutely $\mathbb{Z}$-graded and relatively $H_1(S^3-f(G))$-graded.

\section{Filtered Graph Floer Homology and the $\tau$ Invariant}
\label{ch:Main}

\subsection{Spatial Graphs and the Chain Complex}

In this section, we will define our filtered graph Floer homology chain complex.  It is defined for balanced spatial graphs.

\begin{definition}
	A transverse spatial graph is \textbf{balanced} if there is an equal number of incoming and outgoing edges at each vertex. 
\end{definition}

\begin{figure}
 	\begin{center}
    	\includegraphics[width=.25\textwidth]{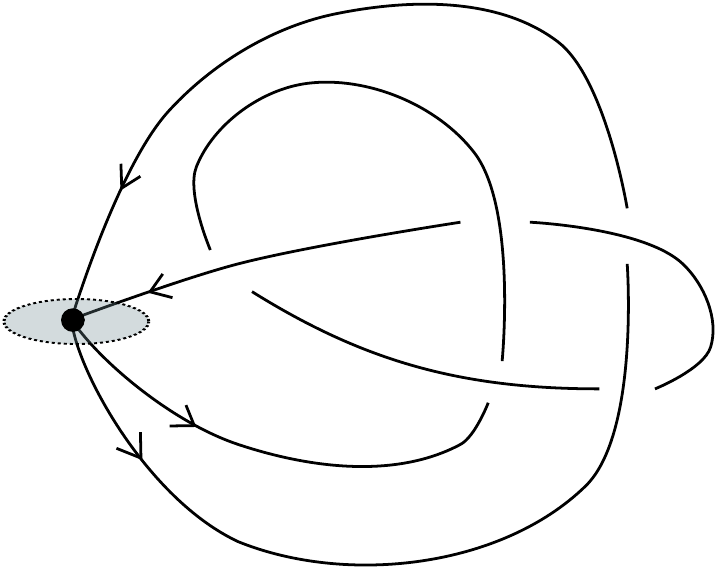}
  	\end{center}
  	\caption{A balanced spatial graph}
	\label{fig:knottedwedge}
\end{figure}

\begin{figure}
	\begin{center}
	\includegraphics[width=.25\textwidth]{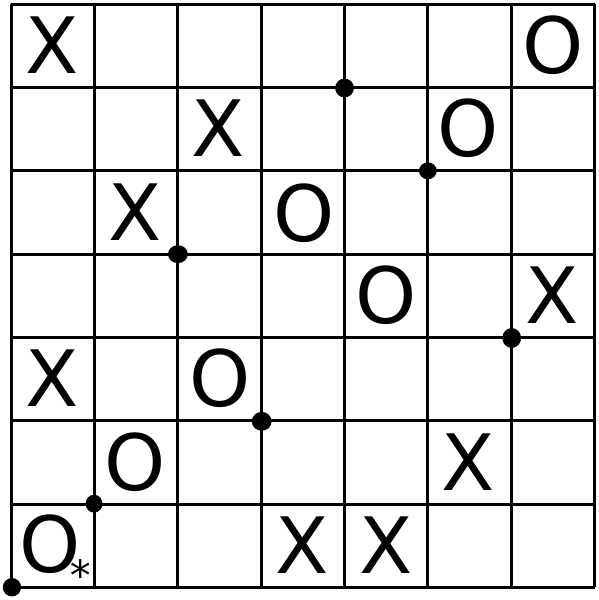}
	\end{center}
	\caption{A generator of $CF^-(g)$}
	\label{fig:gridgen1}
\end{figure}

For an index $n$ grid diagram $g$ representing a spatial graph $f: G \rightarrow S^3$, we choose an ordering for the $O$-markings of $g$ and denote them $O_1, \ldots , O_n$. The chain complex $CF^-(g)$ is freely generated over $\mathbb{F}[U_1, \ldots, U_n]$, where $\mathbb{F} = \mathbb{Z}/2\mathbb{Z}$ and each $U_i$ is a formal variable corresponding to the $O$-marking $O_i$.  It is generated by the set $\mathbf{S}$ of unordered $n$-tuples of intersection points in $g$ with one point on each horizontal and vertical gridline.  The generating set $\mathbf{S}$ is in bijection with $S_n$, the set of permutations of $n$ elements, so $\mathbf{S} = \{ \mathbf{x} = (x_1, \ldots , x_n) | x_i \in \alpha_i \cap \beta_{\sigma(i)} \text{ for some } \sigma \in S_n\}$.  See \cref{fig:gridgen1} for an example of a generator.  

\begin{definition}
 A rectangle $r$ in the grid diagram connects a generator $\mathbf{x}$ to another generator $\mathbf{y}$ if its lower left and upper right corners are points in $\mathbf{x}$, its upper left and lower right corners are points in $\mathbf{y}$, and all other points in $\mathbf{x}$ and $\mathbf{y}$ coincide.  Such a rectangle is empty if its interior does not contain any points of $\mathbf{x}$ and $\mathbf{y}$.  An empty rectangle may contain $X$- and $O$-markings. The set of empty rectangles from $\mathbf{x}$ to $\mathbf{y}$ is  denoted $\mathcal{R}^{\circ} (\mathbf{x},\mathbf{y})$.
\end{definition}

\begin{figure}
	\begin{center}
	\includegraphics[width=.25\textwidth]{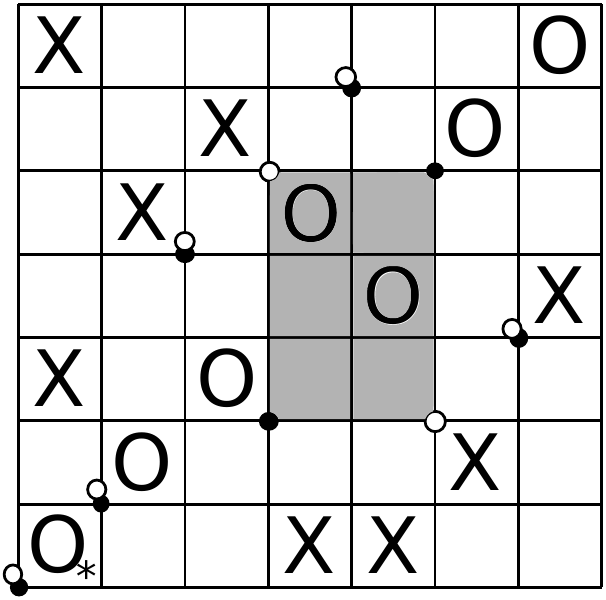}
	\end{center}
	\caption{An empty rectangle connecting the black generator to the white generator}
	\label{fig:gridgen2}
\end{figure}

The boundary map $\partial ^-$ is defined as follows on the generators and extended linearly to $CF^-(g)$:
\[
	\partial^- \mathbf{x} = \sum_{\mathbf{y} \in \mathbf{S}} \mathbf{y} \sum_{r \in \mathcal{R}^{\circ}(\mathbf{x},\mathbf{y})} U_1^{O_1(r)} \cdots U_n^{O_n(r)}
\]
where $O_i (r) = 1$ if $O_i$ is contained in $r$ and $0$ otherwise.

If $g$ is a graph grid diagram representing a balanced spatial graph, the chain complex $CF^-(g)$ is bigraded over $\mathbb{Z}$.  The gradings are defined using the following bilinear map $\mathcal{J}$.

For a point $a = (a_1,a_2)$ and a finite set $B$ of points in the plane, define $\mathcal{J}(a,B)$ to be half of the number of points in $B$ which lie either above and to the right of $a$ or below and to the left of $a$.  That is, $\mathcal{J}(a,B) = \frac 1 2 ( \# \{(b_1,b_2) \in B \, | \text{ either } (a_1 < b_1, a_2 < b_2) \text{ or } (a_1 > b_1, a_2 > b_2) \})$. By extending $\mathcal{J}$ bilinearly to formal sums and differences of sets of points in the plane, we can make the following definition, which is the same as the Maslov grading defined in \cite{MR2372850} and \cite{Harvey-ODonnol}.

\begin{definition}
The Maslov grading, also known as the homological grading, is defined as follows on the generators of the chain complex:
	\[ M(\mathbf{x})=\mathcal{J}(\mathbf{x}-\mathbb{O}, \mathbf{x}-\mathbb{O}) +1\]
where $\mathbb{O}$ and $\mathbb{X}$ are the sets whose points are the $O$- and $X$-markings, respectively.
The Maslov grading is extended to the rest of the chain complex by 
	\begin{align*}
	M(U_i)&=-2 \text{ for all }i \\
	M(0)&=M(1)=0.
	\end{align*}
\end{definition}

For example, the Maslov grading of the element $U_2 U_3^2 \mathbf{x}$ is $M(U_2 U_3^2 \mathbf{x}) = M(\mathbf{x}) - 6$.

\begin{definition}
	The $\mathbb{Z}$-valued Alexander grading is defined as follows for grids which represent balanced spatial graphs (for grids representing spatial graphs that are not balanced, an $H_1(S^3 \setminus f(G))$-valued Alexander grading can be defined, as in \cite{Harvey-ODonnol}):
$$A(\mathbf{x})= \mathcal{J}\left( \mathbf{x} , \mathbb{X}-\sum_{O_i \in \mathbb{O}}m_i O_i \right) $$
where $m_i$ is the weight of $O_i$: the number of $X$-markings in the same column (or equivalently, since we are restricting to balanced graphs, row) as $O_i$.  The Alexander grading is extended to the rest of the chain complex by 
	\begin{align*}
	A(U_i)&=-m_i \text{ for all }i \\
	A(0)&=A(1)=0.
	\end{align*}
	\label{def:Alexgrading}
\end{definition}

We can also view the Alexander grading as a relative grading,  namely $A(\mathbf{x}) - A(\mathbf{y})$, where $\mathbf{x}, \mathbf{y}$ are elements of the chain complex, computed using rectangles.  Any two generators in $\mathbf{S}$ are connected by a sequence of rectangles. This follows from the fact that $\mathbf{S}$ is in bijection with the symmetric group on $n$ letters, $S_n$.  If $\sigma_1, \sigma_2 \in S_n$, there exists a finite sequence of transpositions that will turn $\sigma_1$ into $\sigma_2$.  If $\mathbf{x}_1, \mathbf{x}_2$ are the generators in $\mathbf{S}$ corresponding to $\sigma_1$ and $\sigma_2$, respectively, then that sequence of transpositions corresponds to a sequence of rectangles connecting $\mathbf{x}_1$ to $\mathbf{x}_2$. The following lemma is very similar to Lemma 4.13 in \cite{Harvey-ODonnol}.
	
\begin{lemma} \label{lem:relAlexgrading}
	If $\mathbf{x}, \mathbf{y}$ are generators of the chain complex and $r$ is a rectangle (not necessarily empty) connecting $\mathbf{x}$ to $\mathbf{y}$, then the relative Alexander grading of $\mathbf{x}$ and $\mathbf{y}$ is
\[
A(\mathbf{x})-A(\mathbf{y})=  \left| \mathbb{X} \cap r \right| - \sum_{O_i \in \mathbb{O} \cap r} m_i .
\]
	
	\begin{proof}	
	By \cref{def:Alexgrading},
\begin{align*}
	A(\mathbf{x}) - A(\mathbf{y}) &= \mathcal{J}\left( \mathbf{x} , \mathbb{X}-\sum_{O_i \in \mathbb{O}}m_i O_i \right) - \mathcal{J}\left( \mathbf{y} , \mathbb{X}-\sum_{O_i \in \mathbb{O}}m_i O_i \right) \\
				&= \mathcal{J}\left( \mathbf{x}, \mathbb{X} \right) - \mathcal{J} \left(\mathbf{y},\mathbb{X}\right) - \left(\sum_{O_i \in \mathbb{O}} m_i\left(\mathcal{J}\left(\mathbf{x},O_i\right) - \mathcal{J}\left(\mathbf{y},O_i\right)\right)\right) \\
				&=  \vphantom{\sum_{O_i \in\mathbb{O}}} \frac 1 2 \left(|\mathbb{X} \cap \left( C \cup D \cup r \cup F \cup G \right)| + |\mathbb{X} \cap \left( B \cup C \cup r \cup E \cup F \right)| \right)  \\
					     & \qquad  - \frac 1 2  \left(|\mathbb{X} \cap \left( B \cup C \cup D \cup F \right)| + |\mathbb{X} \cap \left( C \cup E \cup F \cup G \right)| \right)  \\
					     & \qquad - \left( \sum_{O_i \in \mathbb{O}\cap (C \cup D \cup r \cup F \cup G)} \frac{m_i}{2} \right) - \left( \sum_{O_i \in \mathbb{O} \cap (B \cup C \cup r \cup E \cup F)} \frac{m_i}{2} \right) \\
					     & \qquad  + \left( \sum_{O_i \in \mathbb{O}\cap (B \cup C \cup D \cup F)} \frac{m_i}{2} \right) + \left(\sum_{O_i \in \mathbb{O} \cap (C \cup E  \cup F \cup G)} \frac{m_i}{2} \right) \\
				&= |\mathbb{X} \cap r | - \sum_{O_i \in \mathbb{O} \cap r} m_i 
\end{align*}
Where $A,B,C,D,E,F,G,H$ and $r$ are the regions of the grid indicated in \cref{fig:relAlexgrid}.

	\begin{figure}
		\begin{center}
		\def\svgwidth{0.25\textwidth}
   	 	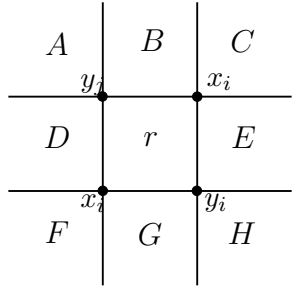
		\end{center}
		\caption{The regions of the grid referred to in \cref{lem:relAlexgrading}}
	\label{fig:relAlexgrid}
	\end{figure}
		\end{proof}
\end{lemma}

\begin{definition}
	The Alexander filtration of $(CF^-(g), \partial ^-)$ is $\{\mathcal{F}^-_m\}_{m \in \mathbb{Z}}$, where $\mathcal{F}^-_m$ is generated by those elements of $CF^-(g)$ whose Alexander grading is less than or equal to $m$.
\end{definition}

\begin{proposition}
	$(CF^-(g), \partial^-)$ is a filtered chain complex.  That is, $\partial^- \circ \partial^- = 0$, the boundary map decreases by one the Maslov grading of elements which are homogeneous with respect to the Maslov grading, and the boundary map preserves the relative Alexander filtration.
	\begin{proof}
		That $\partial^- \circ \partial^- = 0$ follows directly from the proof of Proposition 2.10 of \cite{MR2372850}, since graph grid diagrams differ from link grid diagrams only in the $X$-markings, and the definition of $\partial^-$ does not involve $X$-markings.
		
		The proof that $\partial^-$ decreases Maslov grading by one is also the same as in \cite{MR2372850}. By their Lemma 2.5, if $r$ is an empty rectangle from $\mathbf{x}$ to $\mathbf{y}$, then $M(\mathbf{x}) - M(\mathbf{y}) = 1 - 2 |\mathbb{O} \cap r |$.  Therefore the term in $\partial^- \mathbf{x}$ corresponding to $\mathbf{y}$ will have Maslov grading
		\begin{align*}
		M(U_1^{O_1(r)} \cdots U_n^{O_n(r)}\mathbf{y}) & = M(\mathbf{y}) - \sum_{i = 1}^n 2n_i(r) \\
								& = M(\mathbf{x}) - 1.
		\end{align*}
		
		To show that $\partial^-$ preserves the relative Alexander filtration, note that if a rectangle $r$ connects $\mathbf{x}$ to $\mathbf{y}$, then $A(\mathbf{y}) = A(\mathbf{x}) - \left| \mathbb{X} \cap r \right| + \sum_{O_i \in \mathbb{O} \cap r} m_i $.  Therefore the term in $\partial^- \mathbf{x}$ corresponding to $\mathbf{y}$ will have Alexander grading
		\begin{align*}
		A(U_1^{O_1(r)} \cdots U_n^{O_n(r)}\mathbf{y}) &= A(\mathbf{y}) - \sum_{O_i \in \mathbb{O} \cap r} m_i  \\
									&=  A(\mathbf{x}) - \left| \mathbb{X} \cap r \right| + \sum_{O_i \in \mathbb{O} \cap r} m_i - \sum_{O_i \in \mathbb{O} \cap r} m_i  \\ 
									&= A(\mathbf{x}) - \left| \mathbb{X} \cap r \right| \\
									&\leq A(\mathbf{x}).
		\end{align*}
	\end{proof}
\end{proposition}

\begin{definition}
	Suppose the $O$-markings in $g$ are numbered so that $O_1, \ldots , O_k$ are edge $O$'s and $O_{k+1}, \ldots, O_n$ are vertex $O$'s.  Let $\mathcal{U}$ be the minimal subcomplex of $CF^-(g)$ containing $U_{k+1}CF^-(g) \cup \cdots \cup U_nCF^-(g)$.  Then $(\widehat{CF}(g), \widehat \partial )$ is the filtered chain complex obtained from $(CF^-(g), \partial^-)$ by setting  $\widehat{CF}(g) = CF^-(g) / \mathcal{U}$ and letting $\widehat{\partial}$ be the map on the quotient induced by $\partial^-$.  We consider $\widehat{CF}(g)$ as a vector space over $\mathbb{F}$.   
\end{definition}

We denote by $\widehat{HFG}(g)$ the homology of the associated graded object of $\widehat{CF}(g)$.  It is finitely generated as a vector space over $\mathbb{F}$, since all of the $U_i$'s act trivially on it (\cite{Harvey-ODonnol} Proposition 4.29).

\subsection{Alexander filtration and the $\tau$ invariant}

For a knot $K$, the Alexander filtration of the knot Floer homology chain complex for $K$ is an absolute grading preserved under the maps associated to the commutation and (de)stabilization grid moves.  For balanced spatial graphs, as discussed elsewhere in this chapter, only the relative Alexander filtration of the graph Floer homology chain complex is preserved under the maps associated to the commutation' and (de)stabilization' grid moves.  Therefore, in order to define a $\tau$ invariant for balanced spatial graphs, we need to fix an absolute Alexander grading and filtration of the graph Floer homology chain complex that will be preserved under the maps associated to all of the graph grid moves.

To do this, we show that the homology of the associated graded complex $\oplus_{m} \widehat{\mathcal{F}}_m / \widehat{\mathcal{F}}_{m-1}$ is non-trivial.  To show this, we appeal to the following lemma.

\begin{lemma}
\label{filteredchaincomplexes}
	Let $(C,\partial)$ be a filtered chain complex with filtration $\{\mathcal{F}_s\}$ of $C$ such that $H_*(C) \not = 0$ and $\bigcap\limits_s \mathcal{F}_s = 0$. If for each homological grading $i$, the chain group $C_i$ is finitely generated, then $H_*(\mathcal{F}_s / \mathcal{F}_{s-1}) \not = 0$ for some $s$.  
	
	\begin{proof}
		Since $H_*(C) \not = 0$ and $H_*(C)  = \oplus H_i(C)$, there exists some $i$ for which $H_i(C) \not = 0$.  Therefore there is some non-zero $x \in C_i$ which is homogeneous with respect to the homological grading $i$, with $\partial x = 0$, and whose homology class is nonzero.  We can then choose the minimal filtration level $s$ so that $x \in \mathcal{F}_s$.
		
		Let $\partial_s : \mathcal{F}_s/\mathcal{F}_{s-1} \rightarrow \mathcal{F}_s/\mathcal{F}_{s-1}$.  Then $\partial_s (x + \mathcal{F}_{s-1}) = \partial x + \mathcal{F}_{s-1} = 0 + \mathcal{F}_{s-1}$. If $x + \mathcal{F}_{s-1}$ is not a boundary in the chain complex $(\mathcal{F}_s / \mathcal{F}_{s-1}, \partial_s)$, then $H_*(\mathcal{F}_s / \mathcal{F}_{s-1}) \not = 0$ and we are done.  
		
		If $x + \mathcal{F}_{s-1}$ is a boundary in $(\mathcal{F}_s / \mathcal{F}_{s-1}, \partial_s)$, then there is some $y \in \mathcal{F}_s$ with $x + \mathcal{F}_{s-1} = \partial_s (y + \mathcal{F}_{s-1}) = \partial y + \mathcal{F}_{s-1}$.  Set $z = x - \partial y \in \mathcal{F}_{s-1}$.  Since $x$ is a cycle, $\partial z = \partial (x - \partial y ) = 0$.  Therefore $z \in \mathcal{F}_{s-1}$ is a cycle and since $x$ and $z$ differ by a boundary, $[z] = [x] \not = 0$ in $H_i(C)$.
		
		We can repeat this process, choosing the minimal filtration level $r \leq s-1$ so that $z \in \mathcal{F}_r$, yielding a cycle $z_1 \in \mathcal{F}_{r-1}$ with $[z_1] = [z] = [x] \not = 0$ in $H_i(C)$.  Iterating this process will produce infinitely many representatives of $[x]$, each in different filtration levels.  This contradicts our hypothesis that for each homological grading $i$, the chain group $C_i$ is finitely generated. 
	\end{proof}
\end{lemma}

Note that the grid chain complex $(\widehat{CF}(g), \widehat{\partial})$ satisfies the condition in \cref{filteredchaincomplexes} that for each Maslov grading level $i$, the chain group $\widehat{CF}(g)_i$ is finitely generated.  This is because all elements of $\widehat{CF}(g)$ are of the form $U_1^{a_1} \cdots U_k^{a_k} \mathbf{x}$ for some generator $\mathbf{x}$ and with $a_j \geq 0$ for all $j$, so 
	\[ M(U_1^{a_1} \cdots U_k^{a_k} \mathbf{x}) = M(\mathbf{x}) - 2 \left( \sum_{j = 1} ^k a_j \right) . \]
Since there are finitely many generators, since $M(\mathbf{x})$ is finite, and since there are only finitely many ways to write a given number $i$ as the sum of finitely many positive integers, the condition is satisfied.  

\begin{definition}
\label{def:HFKfilt}
	For a grid diagram $g$ representing a balanced spatial graph $f: G \rightarrow S^3$, define the symmetrized Alexander filtration $\{ \widehat{\mathcal{F}}^H_m\}_{m \in \frac 1 2 \mathbb{Z}}$ to be the absolute Alexander filtration obtained by fixing the relative Alexander grading so that $m_{max}(g) = -m_{min}(g)$, where $m_{max}(g) = max \lbrace m | H_*(\widehat{\mathcal{F}}_{m}(g) / \widehat{\mathcal{F}}_{m-1}(g)) \neq 0 \rbrace$ and $m_{min}(g) = min \lbrace m | H_*(\widehat{\mathcal{F}}_{m}(g) / \widehat{\mathcal{F}}_{m-1}(g)) \neq 0 \rbrace$.
\end{definition}

Now that we have symmetrized the Alexander filtration of $\widehat{CF}(g)$, we can lift that filtration to a symmetrized filtration of $CF^-(g)$.

\begin{definition}
Define the symmetrized Alexander filtration of $CF^-(g)$ to be $\{ \mathcal{F}^{-H}_m\}_{m \in \frac 1 2 \mathbb{Z}}$, obtained by fixing the relative Alexander grading of $CF^-(g)$ so that each generator $\mathbf{x} \in \mathbf{S}(g)$ is in the same filtration level of $\{ \mathcal{F}^{-H}_m\}_{m \in \frac 1 2 \mathbb{Z}}$ as it is in $\{ \widehat{\mathcal{F}}^H_m\}_{m \in \frac 1 2 \mathbb{Z}}$.
\end{definition}

\begin{remark}
This is not necessarily the only way to symmetrize the Alexander filtration.  If we knew that the bigraded Euler characteristic of $\widehat{HFG}(g)$ (which is an Alexander polynomial, see \cite{Harvey-ODonnol}) were non-zero, then we could fix an absolute Alexander grading so that the maximal and minimal terms with non-zero coefficients in the Alexander polynomial were centered around zero.  It would be interesting to answer the question of whether these two ways of fixing the Alexander grading are equivalent.  
\end{remark}

\begin{definition}
\label{def:tau}
For a graph grid diagram $g$ representing a balanced spatial graph $f:G \rightarrow S^3$, define 	the $\tau$ invariant of $g$ to be	\[ \tau (g) = \text{min} \{ m \in \frac 1 2 \mathbb{Z} | \iota _m \text{ is non-trivial} \}\]
	where $\iota_m : H_*(\widehat{\mathcal{F}}_m^H) \rightarrow H_*(\widehat{CF}(g))$ is the map induced by inclusion.
\end{definition}

In proving the next theorem, we will appeal to this lemma.
\begin{lemma}[\cite{rice.160793320010101}, Theorem 3.2]
\label{lem:McCleary}
	If $F: B \rightarrow C$ is a filtered chain map which induces an isomorphism on the homology of the associated graded objects of $B$ and $C$, then $F$ is a filtered quasi-isomorphism.
\end{lemma}

\begin{theorem}
	\label{thm:main}
	For grid diagrams $g, g'$ representing $f:G \rightarrow S^3$,  there exist filtered quasi-isomorphisms  $\phi_1: CF^-(g) \rightarrow CF^-(g')$ and $\phi_2: CF^-(g') \rightarrow CF^-(g)$ which preserve the symmetrized filtration $\lbrace \mathcal{F}^{-H}_m \rbrace$.
	
	\begin{proof}
	For graph grid diagrams $g$ and $g'$ both representing a balanced spatial graph $f: G \rightarrow S^3$, we know by \cref{thm:gridmoves} \cite{Harvey-ODonnol} that there is a finite sequence of cyclic permutation, commutation', stabilization', and destabilization' moves which turns $g$ into $g'$.  Thus, once we show that each of these grid moves is associated to a quasi-isomorphism of filtered chain complexes, we can take the composition of the maps associated to each of the grid moves in the sequence, resulting in a filtered quasi-isomorphism from $CF^-(g)$ to $CF^-(g')$. The proof that each of the grid moves is associated to a quasi-isomorphism of filtered chain complexes consists of three steps:
	
	\begin{enumerate}
		\item We need to show that if $g$ and $\overline{g}$ are graph grid diagrams which are related by a cyclic permutation, commutation', stabilization', or destabilization' grid move, there exists a chain map $\Phi : CF^-(g) \rightarrow CF^-(\overline{g})$ and an integer $s$ such that for all $m$, we have $\Phi (\mathcal{F}^-_m(g)) \subset \mathcal{F}^-_{m+s}(\overline{g})$, and such that $\Phi$ induces an isomorphism $H_*(\mathcal{F}^-_{m}(g) / \mathcal{F}^-_{m-1}(g)) \rightarrow H_*(\mathcal{F}^-_{m+s}(\overline{g}) / \mathcal{F}^-_{m+s-1}(\overline{g}))$.  Note that here, we are working with the original Alexander filtration rather than the symmetrized version.  This will be proved in \cref{sec:cyclicperm}, \cref{sec:comm'}, and \cref{sec:stab'}.
		
		\item We need to show that each of the maps from Step (1) induces a quasi-isomorphism on the symmetrized Alexander filtration.  That is, we need to show that $\Phi_*: H_*(\mathcal{F}_{m}^{-H}(g) / \mathcal{F}_{m-1}^{-H}(g)) \rightarrow H_*(\mathcal{F}_{m}^{-H}(\overline{g}) / \mathcal{F}_{m-1}^{-H}(\overline{g}))$ is an isomorphism.  Since we know from Step (1) that $\Phi$ induces an isomorphism on the homology of the associated graded objects, it is sufficient to show that the span $m_{max} - m_{min}$ is the same for both $\widehat{\mathcal{F}}_m(g)$ and $\widehat{\mathcal{F}}_m(\overline{g})$. We will show that $m_{max}(\overline{g}) = m_{max}(g) + s$ and $m_{min}(\overline{g}) = m_{min}(g) + s$.
		
		Assume for the sake of contradiction that $m_{max}(\overline{g}) > m_{max}(g) + s$.  Then there exists some $\mathbf{y} \in \widehat{\mathcal{F}}_{m_{max}(\overline{g})}(\overline{g})$ such that $[\mathbf{y}] \in H_*(\widehat{\mathcal{F}}_{m_{max}(\overline{g})}(\overline{g})/\widehat{\mathcal{F}}_{m_{max}(\overline{g})-1}(\overline{g}))$ is non-trivial.  Then, since \[\hspace{4em} \Phi_* : H_* (\widehat{\mathcal{F}}_{m_{max}(\overline{g})-s}(g)/\widehat{\mathcal{F}}_{m_{max}(\overline{g})-s-1}(g)) \rightarrow H_* (\widehat{\mathcal{F}}_{m_{max}(\overline{g})}(\overline{g})/\widehat{\mathcal{F}}_{m_{max}(\overline{g})-1}(\overline{g}))\] is an isomorphism, there exists some non-trivial $[\mathbf{x}] = \Phi^{-1}_*([\mathbf{y}])$ in\break  $H_*(\widehat{\mathcal{F}}_{m_{max}(\overline{g})-s}(g)/\widehat{\mathcal{F}}_{m_{max}(\overline{g})-s-1}(g))$.  This contradicts our assumption that $m_{max}(\overline{g}) > m_{max}(g) + s$, so we have that $m_{max}(\overline{g}) \leq m_{max}(g) + s$.  Similar arguments show that $m_{max}(\overline{g}) \geq m_{max}(g) + s$, and that $m_{min}(\overline{g}) = m_{min}(g) + s$.  Therefore we have shown that $m_{max}(g) - m_{min}(g) = m_{max}(\overline{g}) - m_{min}(\overline{g})$.
		
		\item We need to know that the existence of a quasi-isomorphism on the associated graded object of a filtered chain complex implies the existence of a filtered quasi-isomorphism on the filtered chain complex.  This is exactly what \cref{lem:McCleary} \cite{rice.160793320010101} says.
	\end{enumerate}
	\end{proof}
\end{theorem}

\begin{lemma}
	Suppose that there exist filtered quasi-isomorphisms $F : \widehat{CF}(g) \rightarrow \widehat{CF}(\overline{g})$ and $F' : \widehat{CF}(\overline{g}) \rightarrow \widehat{CF}(g)$.  Then $\tau (g) = \tau (\overline{g})$.
	
	\begin{proof}
		Suppose that $\tau(g) = a$.  Then we have the following commutative diagram:
		
			\begin{center}
			\begin{tikzcd}
				H_*(\widehat{\mathcal{F}}_a(g)) \arrow{r}{i_*} \arrow{d}{F^a_*} & H_*(\widehat{CF}(g)) \arrow{d}{F_*} \\
				H_*(\widehat{\mathcal{F}}_a(\overline{g})) \arrow{r}{j_*}  			& H_*(\widehat{CF}(\overline{g})) 
			\end{tikzcd}
			 \end{center}
			 
		Thus there is some $x \in H_*(\widehat{\mathcal{F}}_a(g))$ which maps via $F_* \circ i_*$ to a non-zero element of $H_*(\widehat{CF}(\overline{g}))$ to which $j_*$ sends $F^a_*(x)$.  Therefore $j_*: H_*(\mathcal{F}_a(\widehat{CF}(\overline{g}))) \rightarrow H_*(\widehat{CF}(\overline{g}))$ is non-trivial, so $\tau (\overline{g}) \leq a$.  
		
		The same argument using $F' : \widehat{CF}(\overline{g}) \rightarrow \widehat{CF}(g)$ says that $\tau(g) \leq \tau(\overline{g})$, so putting the two inequalities together gives the result that $\tau(g) = \tau(\overline{g})$.
	\end{proof}
\end{lemma}

With the previous lemma, we have shown the following corollary to \cref{thm:main}.

\begin{corollary}
\label{cor:tau}
	If $g$ and $\overline{g}$ are graph grid diagrams representing a balanced spatial graph $f: G \rightarrow S^3$, then $\tau (g) = \tau (\overline{g})$.
\end{corollary}

Now we have a well-defined $\tau$ invariant for balanced spatial graphs.

\begin{definition}
\label{def:spatialtau}
	For a balanced spatial graph $f: G \rightarrow S^3$, if $g$ is any graph grid diagram representing $f$, then
	\[ \tau (f) = \tau (g). \]
\end{definition}

\subsection{Cyclic Permutation}
\label{sec:cyclicperm}
Suppose that $g$ and $\overline{g}$ are graph grid diagrams which differ by a cyclic permutation move.  Since the chain complex $(CF^-(g), \partial^-_g)$ and $CF^-(\overline{g}), \partial^-_{\overline{g}})$ are defined from toroidal grid diagrams, the chain map associated to the cyclic permutation grid move is the identity map, so it is a quasi-isomorphism.  However, we still need to show that the map preserves the Alexander filtration, which was defined using planar grid diagrams.

From \cref{lem:relAlexgrading} and Corollary 4.14 in \cite{Harvey-ODonnol}, we know that the relative Alexander grading is well-defined on the toroidal grid diagram.  Define new gradings $A'_{g}(\cdot)$ and $A'_{
\overline{g}}(\cdot)$ by shifting the Alexander gradings on $CF^-(g)$ and $CF^-(\overline{g})$, respectively, so that in each one, $\mathbf{x}_{\mathbb{O}}$, the generator whose points are at the lower left corner of each of the grid squares containing and $O$-marking, has grading zero.  Now the identity map preserves this shifted grading.  If $s$ and $\overline{s}$ were the shifts from $A_g(\cdot)$ to $A'_{g}(\cdot)$ and from $A_{\overline{g}}(\cdot)$ to $A'_{\overline{g}}(\cdot)$, respectively, then we see that the identity map sends elements of $CF^-(g)$ with Alexander grading $m$ to elements of $CF^-(\overline{g})$ with Alexander grading $m + s - \overline{s}$.  Therefore $H_*(\mathcal{F}^{-}_{m}(g) / \mathcal{F}^-_{m-1}(g)) \rightarrow H_*(\mathcal{F}^-_{m+s - \overline{s}}(\overline{g}) / \mathcal{F}^-_{m+s - \overline{s}-1}(\overline{g}))$, the map induced by the identity, is an isomorphism.  

\subsection{Commutation'}
\label{sec:comm'}
\begin{figure}
 	\begin{center}
    	\includegraphics[width=.25\textwidth]{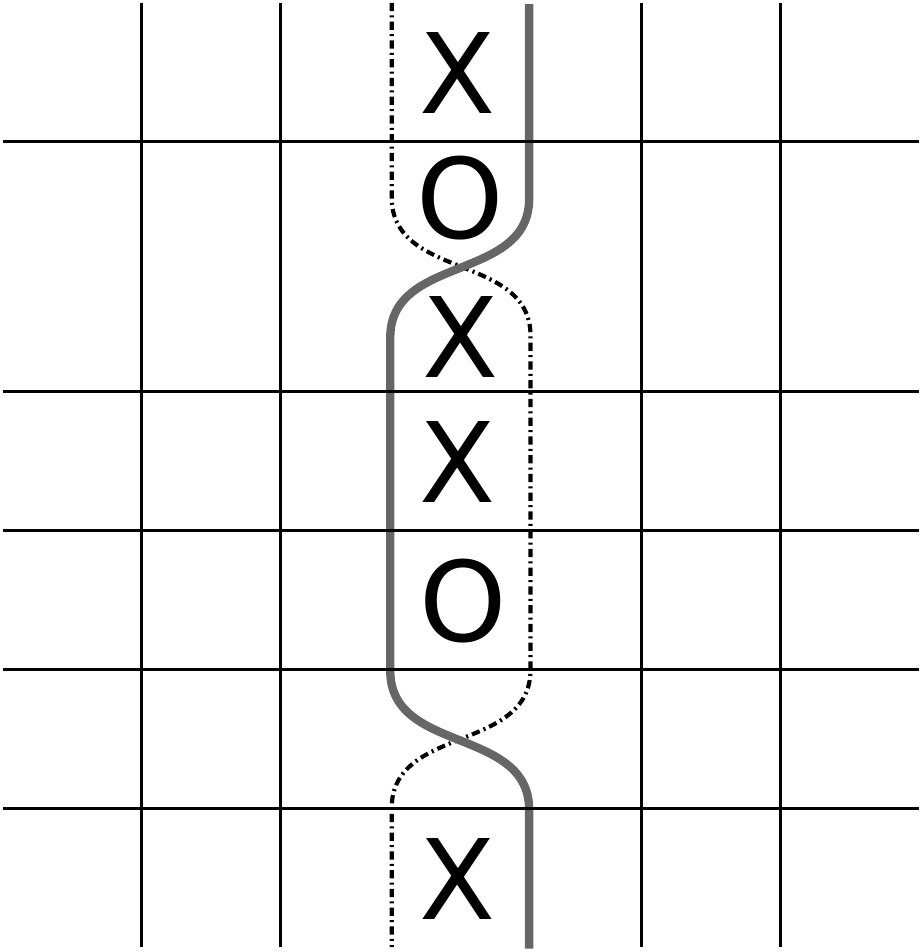}
		\put(-120,-15){$\cdots\: \beta_{n-1}$} \put(-80,-6){$\beta$} \put(-60,-6){$\gamma$} \put(-45,-15){$\beta_1$} \put(-25,-15){$\beta_2$}
  	\end{center}
  	\caption{A grid showing both $g$ and $\overline{g}$.  The grid diagram with $\beta$ but not $\gamma$ is $g$, and the diagram with $\gamma$ but not $\beta$ is $\overline{g}$.}
	\label{fig:CommutationGrid}
\end{figure}
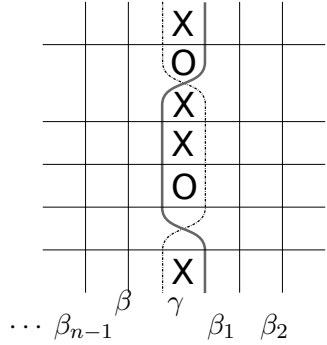

Let $g$ and $\overline{g}$ be graph grid diagrams which differ by a commutation' move.  We can depict both grids in a single diagram, as shown in \cref{fig:CommutationGrid}.  In this example $g$ is the graph grid diagram obtained from \cref{fig:CommutationGrid} by deleting the line labeled $\gamma$, and $\overline{g}$ is the graph grid diagram obtained from it by deleting $\beta$.  The proof of commutation' invariance closely follows that in \cite{Harvey-ODonnol}.

\begin{figure}
 	\begin{center}
    	\includegraphics[width=.25\textwidth]{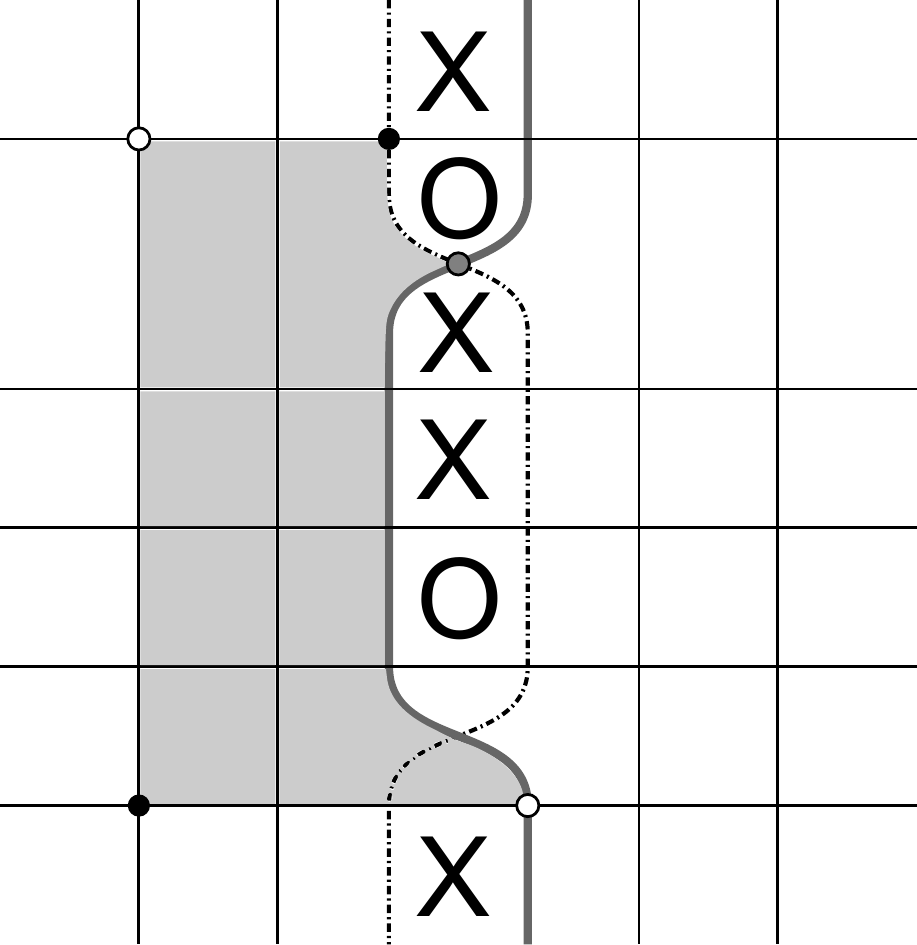}
  	\end{center}
  	\caption{A pentagon connecting the black generator to the white generator, counted in $\Phi'_{\beta \gamma}$.}
	\label{fig:Pentagon}
\end{figure}

Recall that the differential map $\partial ^- : CF^-(g) \rightarrow CF^-(g)$ counts empty rectangles connecting generators in $g$.  In this section, we will consider maps that count empty pentagons and hexagons in the combined grid showing both $g$ and $\overline{g}$.  An embedded pentagon $p$ in the combined grid diagram connects $\mathbf{x} \in \mathbf{S}(g)$ to $\mathbf{y} \in \mathbf{S}(\overline{g})$ if $\mathbf{x}$ and $\mathbf{y}$ agree in all but two points, and if the boundary of $p$ is made up of arcs of five grid lines, whose intersection points are, in counterclockwise order, $a,x_2, y_2, x_1, y_1$, where $a \in \beta \cap \gamma$, $y_1 = \mathbf{y} \cap \gamma$, and $x_1 = \mathbf{x} \cap \beta$.  See \cref{fig:Pentagon} for an example.  Such a pentagon $p$ is \emph{empty} if its interior does not contain any points of $\mathbf{x}$ or $\mathbf{y}$.  The set of empty pentagons connecting $\mathbf{x}$ to $\mathbf{y}$ is denoted $Pent^\circ _{\beta \gamma} (\mathbf{x}, \mathbf{y})$.  

\begin{definition}
For $\mathbf{x} \in S(g)$, let $$\Phi '_{\beta \gamma} (\mathbf{x}) = \sum_{\mathbf{y} \in S(\overline{g})}\left( \sum_{p \in \text{Pent}^\circ_{\beta \gamma}} U_1^{O_1(p)} \dots U_n^{O_n(p)} \cdot \mathbf{y}\right)$$
and note that $\Phi '_{\beta \gamma} (\mathbf{x}) \in CF^-(\overline{g})$.
\end{definition}

\begin{lemma}
The map $\Phi'_{\beta \gamma}$ is a chain map which preserves Maslov grading and respects the Alexander filtration, which is to say that $\Phi'_{\beta \gamma} (\mathcal{F}^-_m(g)) \subset \mathcal{F}^{-}_{m+d}(\overline{g})$ for some $d \in \mathbb{Z}$, where $\{\mathcal{F}^-_m(g)\}$ is the unsymmetrized Alexander filtration of $CF^-(g)$.  Moreover, it induces an isomorphism on the homology of the associated graded object, so 
\[(\Phi'_{\beta \gamma})_*: H_*(\mathcal{F}^-_m(g)/\mathcal{F}^-_{m-1}(g)) \rightarrow H_*(\mathcal{F}^-_{m+d}(\overline{g})/\mathcal{F}^-_{m+d-1}(\overline{g}))\]
is an isomorphism for all $m$.
\begin{proof}
This proof has three parts:
	\begin{enumerate}
		\item $\Phi'_{\beta \gamma}$ preserves Maslov grading.  This follows immediately from Lemma 5.2 in \cite{Harvey-ODonnol} because the difference between their $\Phi_{\beta \gamma}$ map between associated graded chain complexes and our filtered map $\Phi'_{\beta \gamma}$ between filtered chain complexes is that in the filtered setting pentagons may contain $X$-markings, but Maslov grading does not involve the $X$-markings on the grid in any way.

		\item The map $\Phi'_{\beta \gamma}$ preserves the Alexander filtration in the sense given in the statement of the lemma and induces an isomorphism on the homology of the associated graded object.  In the proof of Lemma 5.2 in \cite{Harvey-ODonnol}, Harvey and O'Donnol show that their map $\Phi_{\beta \gamma}$ shifts the Alexander grading by some fixed element $\delta(g, \overline{g}) \in H_1(S^3 \setminus f(G))$, which is the class in $H_1(S^3 \setminus f(G))$ of the sums of the meridians of the graph arcs connecting the $X$ and $O$-markings in the upper region and the lower region between $\beta$ and $\gamma$ in the combined grid.   By collapsing their Alexander grading using the obvious map from $H_1(S^3 \setminus f(G))$ to $\mathbb{Z}$, we obtain from their $\Phi_{\beta \gamma}$ the induced map of our $\Phi'_{\beta \gamma}$ on the associated graded objects $\bigoplus_m \mathcal{F}^-_{m}(g)/\mathcal{F}^-_{m-1}(g) \rightarrow \bigoplus_m \mathcal{F}^-_{m+d}(\overline{g})/\mathcal{F}^-_{m+d-1}(\overline{g})$, where $d \in \mathbb{Z}$ corresponds to $\delta(g, \overline{g}) \in H_1(S^3 \setminus f(G))$ and is the number of graph arcs connecting the $X$ and $O$-markings in the upper region and the lower region between $\beta$ and $\gamma$ in the combined grid.  Note that $d$ may be positive or negative depending on the orientation of the graph arcs.    Therefore we know that $\Phi'_{\beta \gamma}$ induces an isomorphism on the homology of the associated graded object.  
		
		It remains to show that $\Phi'_{\beta \gamma}(\mathcal{F}^-_m(g)) \subset \mathcal{F}^-_{m+d}(\overline{g})$.  Notice that $\Phi'_{\beta \gamma}$ can be decomposed into a sum of $\Phi_{\beta \gamma}$ plus terms corresponding to empty pentagons that contain $X$-markings.  Harvey and O'Donnol use their generalized winding number definition of the Alexander grading to show that $\Phi_{\beta \gamma}$ preserves the Alexander grading.  We need to show that for $\mathbf{x} \in \mathcal{F}^-_m(g)$, each term of $\Phi'_{\beta \gamma}(\mathbf{x})$ corresponding to an empty pentagon containing at least one $X$-marking has Alexander grading less than or equal to $m+d$ in $CF^-(\overline{g})$.  We will also use Harvey and O'Donnol's generalized winding number function, $h(\cdot)$.
		
		As Harvey and O'Donnol did for pentagons not containing any $X$-markings, we will consider $h^g (x_2)-h^g (u_2)$ and $h^{\overline{g}}(y_2) - h^{\overline{g}}(u_1) = h^{\overline{g}}(y_2) - h^g (u_1)$.  For each of these quantities, there are several cases to consider.  However, since we are working with planar, not toroidal, graph grid diagrams, we do not need the cases where $x_2 \in \alpha_i$ for $ i \leq l$ or where $x_1 \in \alpha_i$ for $i > k+1$.
		
\begin{figure}
\def\svgwidth{.25\textwidth}
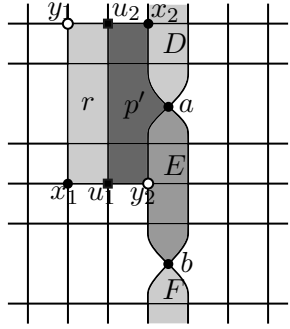
\caption{A pentagon $p$ composed of a rectangle $r$ and a narrow pentagon $p'$.}
\label{fig:PentagonRegions}
\end{figure}		
		
		First, we look at $h^g (x_2)-h^g (u_2)$ when $x_2 \in \alpha_i$ for $i > k$.  In this case, there are no graph arcs passing between $x_2$ and $u_2$, so $h^g (x_2)-h^g (u_2)=0$.  We note that the narrow pentagon $p'$ has empty intersection with the region marked $E$ in \cref{fig:PentagonRegions}, so $h^g (x_2)-h^g (u_2)=0 = - \left| \mathbb{X} \cap(p' \cap E) \right| + \sum_{O_i \in p' \cap E} m_i$.  In the case that $x_2 \in \alpha_i$ for $l+1 \leq i < k$, there may be graph arcs passing between $x_2$ and $u_2$.  If the $O$-marking in $E$ is in $p'$, then there will be one downward-pointing graph arc for each $X$-marking in $E \setminus p'$, and if the $O$-marking is in $E \setminus p'$, then there will be one upward-pointing graph arc for each $X$-marking in $E \cap p'$.  In either case, we see that $h^g (x_2)-h^g (u_2) = - \left| \mathbb{X} \cap(p' \cap E) \right| + \sum_{O_i \in p' \cap E} m_i$.
		
		Now we consider $h^{\overline{g}}(y_2) - h^{\overline{g}}(u_1)$.  Let $d$ be the number of graph arcs from region $D$ to region $F$ as marked in \cref{fig:PentagonRegions}.  It is negative if those arcs are downward-pointing and positive if they are upward-pointing.  Note that $d$ is the image of  $\delta(g, \overline{g}) \in H_1(S^3 \setminus f(G))$ from the proof of Lemma 5.2 in \cite{Harvey-ODonnol} under the map sending the meridian of a graph edge in $H_1(S^3 \setminus f(G))$ to $1 \in \mathbb{Z}$.  If $x_1 \in \alpha_i$ for $i \leq l$, there are three possibilities for the location of the $O$-marking $O''$ in column $n-1$ of $\overline{g}$: $O'' \in p'$, $O'' \in D$, or $O'' \in F \setminus p'$.  
		
		If $O'' \in p'$, then there is one upward-pointing graph arc between $u_1$ and $y_2$ for each $X$-marking in $F \setminus p'$.  The number of $X$-markings in $F \setminus p'$ is $m'' - d - \left| \mathbb{X} \cap (p' \cap F) \right|$, where $m''$ is the multiplicity of $O''$ and we note that $d = \left| \mathbb{X} \cap D \right|$.  So $h^{\overline{g}}(y_2) - h^{\overline{g}}(u_1) = - \sum_{O_i \in p' \cap F} m_i + d + \left| \mathbb{X} \cap (p' \cap F) \right|$.  If $O'' \in D$, then as in the previous case there is one upward-pointing graph arc between $u_1$ and $y_2$ for each $X$-marking in $F \setminus p'$.  There are $m'' - \left| \mathbb{X} \cap D \right| - \left| \mathbb{X} \cap (p' \cap F) \right|$ such markings, and in this case we note that $d = -m'' + \left| \mathbb{X} \cap D \right|$, so the number of upward-pointing graph arcs between $u_1$ and $y_2$ is $-d - \left| \mathbb{X} \cap (p' \cap F) \right| + \sum_{O_i \in p' \cap F} m_i$, where the sum is empty.  Therefore $h^{\overline{g}}(y_2) - h^{\overline{g}}(u_1) = - \sum_{O_i \in p' \cap F} m_i + d + \left| \mathbb{X} \cap (p' \cap F) \right|$.  If $O'' \in F \setminus p'$, then there is one downward-pointing graph arc between $u_1$ and $y_2$ for each $X$-marking in $D \cup (p' \cap F)$.  We notice that in this case $d = \left| \mathbb{X} \cap D \right|$, so the number of downward-pointing graph arcs between $u_1$ and $y_2$ is $d + \left| \mathbb{X} \cap (p' \cap F) \right| - \sum_{O_i \in p' \cap F} m_i$, where the sum is empty.  So $h^{\overline{g}}(y_2) - h^{\overline{g}}(u_1) = - \sum_{O_i \in p' \cap F} m_i + d + \left| \mathbb{X} \cap (p' \cap F) \right|$.
		
		We see from the above that in all cases, 
		\[ h^g (x_2)-h^g (u_2) = - \left| \mathbb{X} \cap(p' \cap E) \right| + \sum_{O_i \in p' \cap E} m_i \] and \[
		h^{\overline{g}}(y_2) - h^{\overline{g}}(u_1) = \left| \mathbb{X} \cap (p' \cap F) \right| - \sum_{O_i \in p' \cap F} m_i + d. \]
		
		We now put these together to consider 
		\begin{align*}
		h^g (x_1) &+ h^g (x_2) - h^{\overline{g}}(y_1) - h^{\overline{g}}(y_2)  \\
		&= \left[ h^g(x_1) + h^g(u_2) - h^g (y_1) - h^g (u_1)\right] \\
		& \qquad						+ \left[ h^g (x_2)-h^g (u_2) - h^{\overline{g}}(y_2) + h^{\overline{g}}(u_1) \right] \\
		&= \left[ - \left| \mathbb{X} \cap r \right| + \sum_{O_i \in r} m_i \right] \\
		&	\qquad		+ \left[ - \left| \mathbb{X} \cap(p' \cap E) \right| + \sum_{O_i \in p' \cap E} m_i  \right. \\
		& \qquad \left. - \left| \mathbb{X} \cap (p' \cap F) \right| + \sum_{O_i \in p' \cap F} m_i - d \right]\\
		&= - \left| \mathbb{X} \cap p \right| + \sum_{O_i \in p} m_i - d.
		\end{align*}
		Since $A^g(\mathbf{x}) = \sum_{x_i \in \mathbf{x}} -h^g(x_i)$, we can see that 
		\[ A^g(\mathbf{x}) - A^{\overline{g}} (\mathbf{y}) = \left| \mathbb{X} \cap p \right| - \sum_{O_i \in p} m_i + d.\]

		\item The map is a chain map, that is $\partial ^- \circ \Phi'_{\beta \gamma} + \Phi'_{\beta \gamma} \circ \partial ^- = 0$.  This follows immediately from the proof of Lemma 3.1 in \cite{MR2372850}.
	\end{enumerate}
\end{proof}
\end{lemma}

\begin{figure}
 	\begin{center}
    	\includegraphics[width=.25\textwidth]{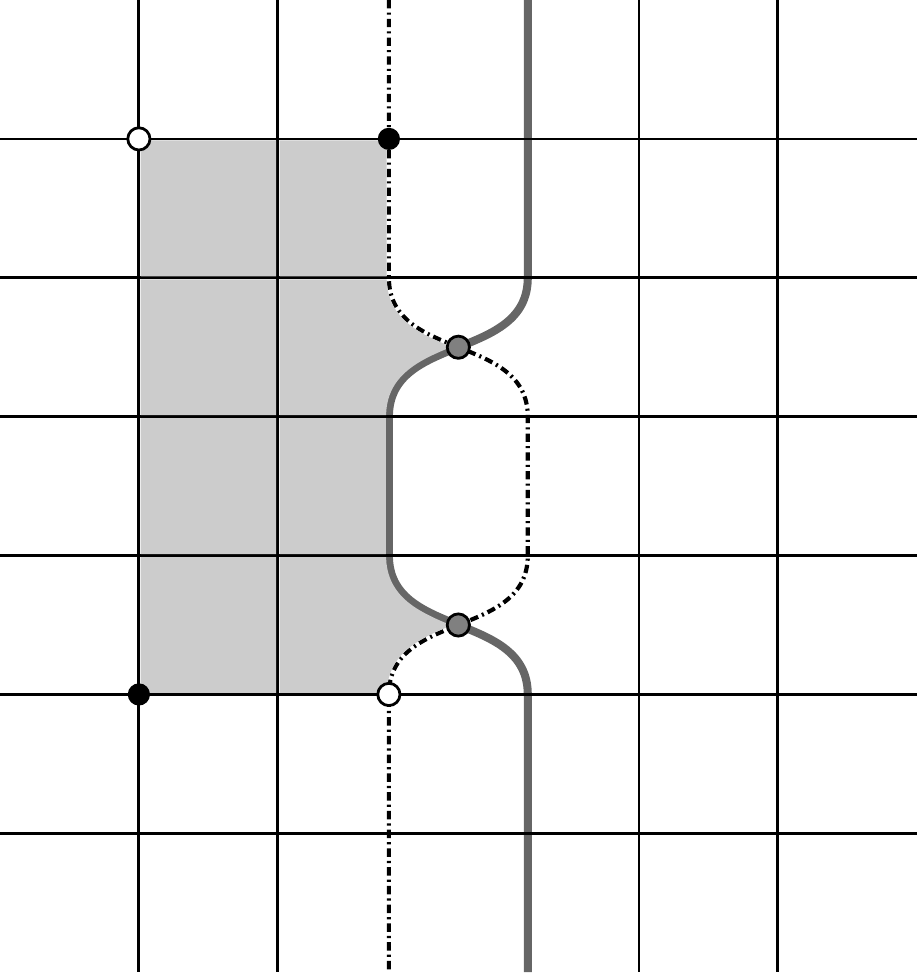}
  	\end{center}
  	\caption{A hexagon counted in $H_{\beta \gamma \beta}$.}
	\label{fig:Hexagon}
\end{figure}

The proof that $\Phi'_{\beta \gamma}$ is a chain homotopy equivalence is the same as the proof in Section 3.1 of \cite{MR2372850}.  An embedded hexagon $h$ in the combined grid showing both $g$ and $\overline{g}$ connects $\mathbf{x} \in \mathbf{S}(g)$ to $\mathbf{y} \in \mathbf{S}(\overline{g})$ if $\mathbf{x}$ and $\mathbf{y}$ agree in all but two points (without loss of generality, say the points where they do not agree are $x_1, x_2$ and $y_1, y_2$), and if the boundary of $h$ is made up of arcs of grid lines whose intersection points are, in counterclockwise order,  $x_1, y_1, a_1, a_2, x_2,$ and $y_2$, where $\lbrace a_1, a_2 \rbrace = \beta \cap \gamma$, and if the interior angles of $h$ are all less than straight angles.  See \cref{fig:Hexagon} for an example.  A hexagon is empty if its interior does not contain any points of $\mathbf{x}$ or $\mathbf{y}$.  The set of empty hexagons connecting $\mathbf{x}$ to $\mathbf{y}$ is denoted $Hex^{\circ}_{\beta \gamma \beta}(\mathbf{x}, \mathbf{y})$.   The chain homotopy operator $H_{\beta \gamma \beta}: CF^-(g) \rightarrow CF^-(g)$ is defined as follows:
	\[H_{\beta \gamma \beta} (\mathbf{x}) = \sum_{\mathbf{y} \in S(g)}\left( \sum_{h \in \text{Hex}^\circ_{\beta \gamma \beta}(\mathbf{x}, \mathbf{y})} U_1^{O_1(h)} \dots U_n^{O_n(h)} \cdot \mathbf{y}\right). \]

\begin{lemma}[\cite{MR2372850} Proposition 3.2]
	The map $\Phi'_{\beta \gamma}$ is a chain homotopy equivalence.  That is,
		\[ \mathbb{I}_{C(\overline{g})} + \Phi'_{\beta \gamma} \circ \Phi'_{\gamma \beta} + \partial^- \circ H_{\gamma \beta \gamma} + H_{\gamma \beta \gamma} \circ \partial^-  = 0 \]
	and 
		\[ \mathbb{I}_{C({g})} + \Phi'_{\gamma \beta} \circ \Phi'_{\beta \gamma} + \partial^- \circ H_{\beta \gamma \beta} + H_{\beta \gamma \beta} \circ \partial^- = 0. \]
		
\end{lemma}

\subsection{Stabilization'}
\label{sec:stab'}
Let $g$ and $g'$ be two graph grid diagrams such that a stabilization' move on $g$ results in $g'$.  Our proof that the (de)stabilization' move induces filtered quasi-isomorphisms in both directions between $CF^- (g)$ and $CF^-(g')$ is modeled on Sarkar's proof in \cite{MR2915478}.  

Sarkar \cite{MR2915478} distinguishes between two types of (de)stabilizations: those at ordinary $O$-markings, which he refers to as $S^3$-grid move (4), and those at special $O$-markings, which he refers to as $S^3$-grid move (5) (Special $O$-markings in the spatial graph case are the vertex $O$'s).  The first type can correspond to a (de)stabilization (Link-grid move (3)), which preserves isotopy class, or a birth in the cobordism (Link-grid move (4)), while the second type corresponds to a death in the cobordism (Link-grid move (7)).

Therefore, although both types of stabilization will be needed to prove the link cobordism result in \cref{thm:LinkCob}, for the purposes of proving the invariance of $HFG^-$ and the $\tau$ invariant, we will only need the first type.

Sarkar defines two stabilization maps, $s_{11}, s_{22}$, and two destabilization maps, $d_{11},d_{22}$.  The $11$ maps correspond to the stabilization in which the new $O$-marking is placed in the row above the $X$ being stabilized, and the $22$ maps correspond to the stabilization in which the new $O$-marking is placed in the row below the $X$ being stabilized.  Because we can use the commutation' move, we only need the graph grid diagram analogs of the $11$ maps. The case for which \cite{MR2915478} uses the $s_{22}, d_{22}$ maps can instead be addressed in the spatial graph case using a commutation' move, then $d_{11}$ or $s_{11}$, then another commutation' move. 

\begin{figure}
			\begin{center}
			\includegraphics{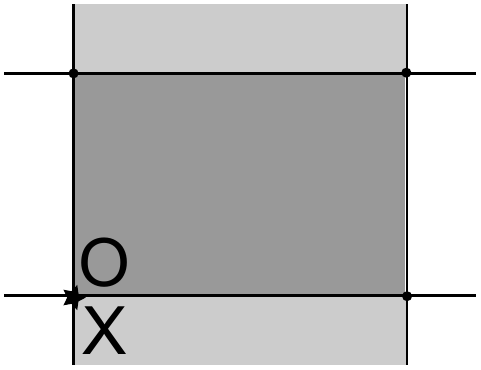}
			\end{center}
			\caption{the rectangles $D$ and $D'$}
			\label{fig:dsmaps}
\end{figure}

The maps $d_{11}: CF^-(g') \rightarrow CF^-(g)$ and $s_{11}:CF^-(g) \rightarrow CF^-(g')$ are defined as follows on the generators of the chain complexes:

\begin{eqnarray*}
	d_{11}(U_{0}^m \mathbf{x}) = U_j^m \sum_{\mathbf{y}} \mathbf{y} \sum_{D \in \mathcal{S}_{1}(\mathbf{x}, \mathbf{y} \cup \star, \star)} U_1^{O_1(D)} \cdots U_n^{O_n(D)} \\
	s_{11}(\mathbf{x}) = \sum_{\mathbf{y}} \mathbf{y} \sum_{D \in \mathcal{S}_{3}(\mathbf{x} \cup \star ,\mathbf{y}, \star)} U_1^{O_1(D)} \cdots U_n^{O_n(D)}
\end{eqnarray*}
 Here, $O_0$ is the the $O$-marking in $g'$ but not in $g$, $O_j$ is the $O$-marking in the row immediately below $O_0$,  $\mathcal{S}_{1}(\mathbf{x}, \mathbf{y} \cup \star, \star)$ and $\mathcal{S}_{3}(\mathbf{x} \cup \star ,\mathbf{y}, \star)$ are the sets of snail-like domains illustrated in Figure 5 of \cite{MR2915478}, and $\star$ is the intersection point of the $\alpha$ and $\beta$ curves immediately below and to the left of the new $O$-marking (see \cref{fig:dsmaps}).

The map $d_{11}$ is exactly the map $F^R$ defined in \cite{MR2372850} and used in \cite{Harvey-ODonnol}, considered as a map from $C$ to $B$, where $C$ is the chain complex associated to the stabilized grid diagram and $B$ is the chain complex associated to the unstabilized grid diagram.  Therefore by Lemma 3.5 in \cite{MR2372850}, the map $d_{11}$ is a chain map which preserves the Maslov grading. In Lemma 5.8 in \cite{Harvey-ODonnol}, Harvey and O'Donnol prove that $d_{11}$ induces an isomorphism $\mathcal{F}^-_{m}(g')/\mathcal{F}^-_{m-1}(g') \rightarrow \mathcal{F}^-_{m+a}(g)/\mathcal{F}^-_{m+a-1}(g)$ for all $m \in \mathbb{Z}$.  When mapped to the integers, the grading shift is $a= -A^{g'}(\star) -1$.  The proof that $d_{11}$ preserves the Alexander filtration up to a shift by $a$ is similar to the proof in \cite{Harvey-ODonnol} that it induces an isomorphism on the associated graded object, except that in the filtered case, we allow the domains to contain $X$-markings, which lowers the Alexander grading of the terms associated to the domains containing $X$-markings. 

\begin{lemma}
	The composition $d_{11} \circ s_{11}$ is the identity map on the associated graded chain complex for the unstabilized grid diagram.
	\begin{proof}
		In the associated graded chain complexes, the only regions counted in $s_{11}$ and $d_{11}$  are rectangles with the starred grid intersection point as their lower left and upper left corners, respectively.  All higher complexity snail-like regions counted in these maps contain the $X$ being stabilized and thus are not counted in the associated graded version.  Furthermore, in the associated graded chain complexes the regions counted may not contain any $X$-markings other than the one in the newly added column.
		
		If $D$ is a rectangle connecting $\mathbf{x} \cup \star$ to $\mathbf{y}$ which is counted in $s_{11}(\mathbf{x})$, then we consider $d_{11}(\mathbf{y})$.  If $D'$ is a domain counted in $d_{11}(\mathbf{y})$, then the boundary of $\partial D' \cap \beta_1$ is $\mathbf{y} - \star$.  Therefore the upper boundary of $D'$ is $\alpha_1$, so the term in $d_{11}(\mathbf{y})$ corresponding to $D'$ is $\mathbf{x}$.  See \cref{fig:dsmaps}. 
		
		No $U_i$'s survive in $d_{11} \circ s_{11}(\mathbf{x})$ since the composite map counts domains $D \cup D'$, which as just discussed are the union of entire columns in the grid diagram.  Since every column contains at least one $X$-marking, the only $D \cup D'$ that may be counted is the single column containing the new $O$-marking.  Since the new $O$-marking is not counted, $d_{11} \circ s_{11}(\mathbf{x})=\mathbf{x}$ and so the composition $d_{11} \circ s_{11}$ is the identity map.
	\end{proof}
\end{lemma}

\begin{lemma}
	The map $s_{11}$ is a quasi-isomorphism between the associated graded chain complexes for the unstabilized and stabilized grid diagrams.
	\begin{proof}
		We know from the previous lemma that $d_{11} \circ s_{11}$ is the identity map on the associated graded chain complex for the unstabilized grid diagram.  The identity map is a quasi-isomorphism, and by Proposition 5.13 in \cite{Harvey-ODonnol}, $d_{11} = F^R$ is a quasi-isomorphism.  Then since $s_{11}$ is the one-sided inverse of a quasi-isomorphism, it is also a quasi-isomorphism.
	\end{proof}
\end{lemma}

\begin{lemma}
	The map $s_{11}$ is a filtered chain map which preserves Maslov grading and respects the Alexander filtration up to a finite shift, so that $s_{11}(\mathcal{F}^-_{m}(g)) \subset \mathcal{F}^-_{m-a}(g')$ and it induces an isomorphism $\mathcal{F}^-_{m}(g)/\mathcal{F}^-_{m-1}(g) \rightarrow \mathcal{F}^-_{m-a}(g')/\mathcal{F}^-_{m-a-1}(g')$ for all $m \in \mathbb{Z}$ and $a= -A^{g'}(\star) -1$.
		\begin{proof}
		By definition, $s_{11}$ is a module homomorphism.  We need to show that it preserves the Maslov grading, it respects the Alexander filtration, and that it is a chain map.
		
		The proof that $s_{11}$ is a chain map is the same as the proof of Lemma 3.5 in \cite{MR2372850} except that the snail-like domains are rotated 90\degree counterclockwise.   
		
		To show that the Maslov grading is preserved, suppose that there is some snail-like domain $D$ which connects $\mathbf{x} \cup \star$ to $\mathbf{z}$ in the stabilized grid $g'$ which is counted in $s_{11}(\mathbf{x})$. We begin by using the definition of the Maslov grading to compare the grading of $\mathbf{x} \cup \star$ in the stabilized grid $g'$ to the grading of $\mathbf{x}$ in the unstabilized grid $g$.
		\begin{eqnarray*}
			M(\mathbf{x} \cup \star) &=& \mathcal{J} (\mathbf{x} + \star - \mathbb{O}_{g'}, \mathbf{x} + \star - \mathbb{O}_{g'} ) + 1 \\
							&=& \mathcal{J}(\mathbf{x}, \mathbf{x}) + \mathcal{J}(\star, \star) + \mathcal{J}(\mathbb{O}_{g'}, \mathbb{O}_{g'}) + 2\mathcal{J}(\mathbf{x}, \star) \\
							&&\hspace{2cm}- 2\mathcal{J}(\mathbf{x}, \mathbb{O}_{g'}) - 2\mathcal{J}( \star, \mathbb{O}_{g'}) +1
		\end{eqnarray*}
		Noting that $\mathbb{O}_{g'}$, the set of $O$-markings in $g'$, is the same as $\mathbb{O}_g \cup O_{k+1}$, where $\mathbb{O}_g$ is the set of $O$-markings in $g$ and $O_{k+1}$ is the new $O$-marking, we can see that 
		\begin{eqnarray*}
			M(\mathbf{x} \cup \star) = \mathcal{J}(\mathbf{x}, \mathbf{x}) &+& \mathcal{J}(\star, \star) + \mathcal{J}(\mathbb{O}_{g}, \mathbb{O}_{g}) + \mathcal{J}(O_{k+1}, O_{k+1}) + 2\mathcal{J}(\mathbb{O}_g, O_{k+1}) \\
					&+& 2\mathcal{J}(\mathbf{x}, \star) - 2\mathcal{J}(\mathbf{x}, \mathbb{O}_{g}) -2\mathcal{J}(\mathbf{x}, O_{k+1}) - 2\mathcal{J}( \star, \mathbb{O}_{g}) \\
					&-& 2\mathcal{J}(\star, O_{k+1}) +1.
		\end{eqnarray*}
		We can use the following observations to simplify the expression:
		\begin{eqnarray*}
			\mathcal{J}(\star, \star) &=& \mathcal{J}(O_{k+1}, O_{k+1}) = 0 \\
			\mathcal{J}(\mathbf{x}, \star ) &=& \mathcal{J}(\mathbf{x}, O_{k+1}) \\
			\mathcal{J}(\mathbb{O}_g, O_{k+1}) &=& \mathcal{J}( \star, \mathbb{O}_{g}) \\
			M(\mathbf{x}) &=& \mathcal{J}(\mathbf{x}, \mathbf{x}) + \mathcal{J}(\mathbb{O}_{g}, \mathbb{O}_{g}) - 2\mathcal{J}(\mathbf{x}, \mathbb{O}_{g}) + 1 \\
			\mathcal{J}( \star, O_{k+1}) &=& \frac 1 2
		\end{eqnarray*}
		Therefore $M(\mathbf{x} \cup \star) = M(\mathbf{x}) -1$.
		
		Since $D$ connects $\mathbf{x} \cup \star$ to $\mathbf{z}$, the term corresponding to $D$ in $s_{11}(\mathbf{x})$ is\break $U_1^{O_1(D)} \cdots U_n^{O_n(D)}\mathbf{z}$.  Therefore we need to compare $M(\mathbf{x})$ and $M(U_1^{O_1(D)} \cdots U_n^{O_n(D)}\mathbf{z})$.  By the definition of Maslov grading, $M(U_1^{O_1(D)} \cdots U_n^{O_n(D)}\mathbf{z}) = M(\mathbf{z}) - 2 \sum_{i=1}^n O_i(D)$, where the summation does not include $i=0$, which corresponds to the new $O$-marking.  Using Lemma 2.5 in \cite{MR2372850}, we know that $M(\mathbf{x} \cup \star) = M(\mathbf{z}) + 1 +2m - 2\sum_{i=0}^{n} O_i(D)$, with $i=0$ included in the sum and where $m$ is the multiplicity of $\star$ in $D$.  Since the multiplicity of $\star$ in the interior of $D$ is $O_{0}(D) - 1$, we can put all of this together to see that
			\[ M(\mathbf{x}) = M(U_1^{O_1(D)} \cdots U_n^{O_n(D)}\mathbf{z}) + 2 - 2 O_{0}(D),\]
so $s_{11}$ preserves the Maslov grading.

		For the Alexander filtration, we need to show that $s_{11}(\mathcal{F}^-_{m}(g)) \subset \mathcal{F}^-_{m-a}(g')$.  Suppose that $D$ is a snail-like domain counted in $s_{11}(\mathbf{x})$.  Then considered in the stabilized grid $g'$, $D$ is a domain connecting $\mathbf{x} \cup \star$ to some generator $\mathbf{z}$.  By Lemma 4.13 in \cite{Harvey-ODonnol}, 
		\[
		A^{g'}(\mathbf{x} \cup \star) - A^{g'}(\mathbf{z}) = n_{\mathbb{X}}(D) - \sum_{i=0}^n m_i \cdot O_{i}(D),
		\]
where $n_{\mathbb{X}}$ is the number of $X$-markings contained in $D$, counted with multiplicity.
		The term of $s_{11}(\mathbf{x})$ corresponding to the domain $D$ is $U_1^{O_1(D)} \cdots U_n^{O_n(D)}\mathbf{z}$, which has Alexander grading $A^{g'}(\mathbf{z}) - \sum_{i = 1}^n O_i(D)$. The shift in the Alexander grading from the $s_{11}$ map is 
\begin{align*}
A^{g'}(U_1^{O_1(D)} \cdots U_n^{O_n(D)}\mathbf{z}) &= A^{g'}(\mathbf{z}) - \sum_{i=1}^n m_iO_i(D) - A^g(\mathbf{x})\\
		&= -n_{\mathbb{X}}(D) + O_0(D) + A^{g'}(\mathbf{x} \cup \star) - A^g(\mathbf{x}) \\
		&= -n_{\mathbb{X}}(D) + O_0(D) + A^{g'}(\mathbf{x}) + A^{g'}(\star) - A^g(\mathbf{x}) \\
		&= -n_{\mathbb{X}}(D) + 1 + A^{g'}(\star).
\end{align*}		
Notice that for domains that do not contain any $X$-markings, which are exactly the domains considered in the associated graded object, the shift in Alexander grading is $1 +A^{g'}(\star)$, which is the negative of the shift for the $d_{11}$ map.  For domains that do contain $X$-markings, the Alexander grading in the terms of $s_{11}(\mathbf{x})$ (for $\mathbf{x}$ in Alexander grading $m$) corresponding to those domains have Alexander grading less than $m +A^{g'}(\star) +1$, since the presence of $X$-markings in the domain reduces their Alexander grading. Therefore $s_{11}(\mathcal{F}^-_m(g)) \subset \mathcal{F}^-_{m-a}(g')$, for $a = -A^{g'}(\star) -1$, and since $d_{11}\circ s_{11}= id$ induces an isomorphism on $\mathcal{F}^-_{m}(g)/\mathcal{F}^-_{m-1}(g) \rightarrow \mathcal{F}^-_{m}(g)/\mathcal{F}^-_{m-1}(g)$ for all $m \in \mathbb{Z}$, we know that $s_{11}$ induces an isomorphism $\mathcal{F}^-_{m}(g)/\mathcal{F}^-_{m-1}(g) \rightarrow \mathcal{F}^-_{m-a}(g')/\mathcal{F}^-_{m-a-1}(g')$ for all $m \in \mathbb{Z}$.

		\end{proof}
\end{lemma}

Using \cref{lem:McCleary} \cite{rice.160793320010101} and the results above that $d_{11}$ and $s_{11}$ are filtered chain maps which are quasi-isomorphisms on the associated graded objects, we see that they are quasi-isomorphisms on the filtered chain complexes on which they are defined.

\section{Link Cobordisms}
\label{ch:Links}

In this section we state the definition of link cobordism and prove an inequality for links analogous to the one proven for knots by Sarkar in \cite{MR2915478}.  This gives an obstruction to sliceness for some links.  In recent independent work, Cavallo \cite{Cavallo} has defined a $\tau$ invariant for links and proven a result similar to \cref{thm:LinkCob}.

\begin{definition}
	A cobordism from a link $L_0$ to another line $L_1$ is a surface $F$ properly embedded in $S^3 \times [0,1]$, such that $F \cap S^3 \times \{0\} = L_0$ and $F \cap S^3 \times \{1\} = - L_1$.  If such a surface exists, we say that $L_0$ is cobordant to $L_1$.  
	
	If two $l$-component links $L_0$ and $L_1$ are connected by a cobordism consisting of $l$ disjoint annuli, we say that they are concordant, and a link which is concordant to the unlink is slice.
\end{definition}

Following \cite{MR2915478}, for the purposes of this section we will allow an $X$-marking and an $O$-marking to occupy the same grid square.  In this case, those two markings represent a small, unknotted link component.  In addition, we call a link grid diagram \emph{tight} if there is exactly one special $O$-marking on each link component.

Also following \cite{MR2915478}, a cobordism between two links can be represented by a series of link grid moves.  These moves are commutations and stabilizations (which correspond to isotopy of links) and births, deaths, $X$-saddles, and $O$-saddles.

\begin{figure}
	\begin{center}
		\subfigure{\includegraphics[width=.25\textwidth]{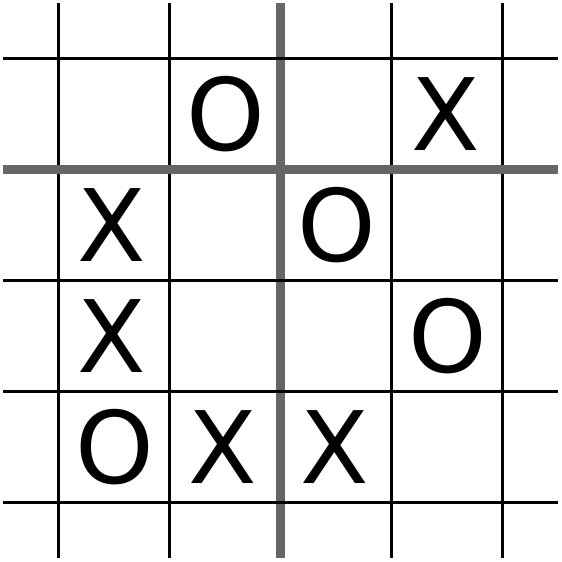}}
		\subfigure{\includegraphics[width=.25\textwidth]{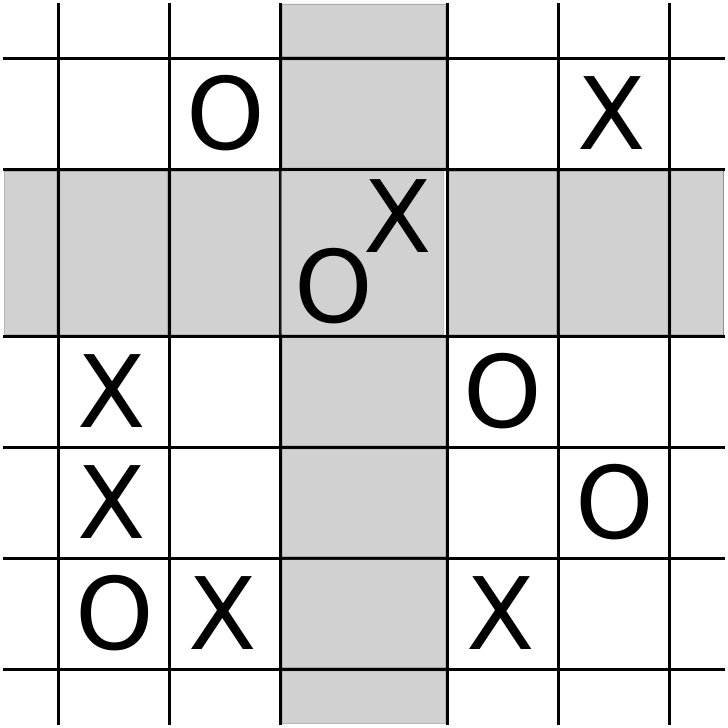}}
	\end{center}
	\caption{A birth in the grid on the left produces the grid on the right; a death in the right-hand grid produces the left-hand grid}
	\label{fig:birth-death}
\end{figure}

A grid diagram $\overline{g}$ is obtained from another grid diagram $g$ via a \textbf{birth} if adding an additional row and column to $g$ and placing both an $O$- and an $X$-marking in the grid square that is the intersection of the new row and column results in $\overline{g}$.  See \cref{fig:birth-death} for an example. This move is link-grid move (4) in \cite{MR2915478}.

A grid diagram $\overline{g}$ is obtained from another grid diagram $g$ via a \textbf{death} if there are a row and a column in $g$, each containing exactly one $X$-marking, whose intersection contains that $X$- and an $O$-marking.  Then $\overline{g}$ is the result of deleting those markings and deformation retracting the row and column to an $\alpha$ and a $\beta$ circle.  For an example, see \cref{fig:birth-death}.  This is very similar to link-grid move (7) in \cite{MR2915478}, with the difference being that Sarkar required the $O$-marking in the dying component to be a special $O$, and here it is a regular $O$-marking.

\begin{figure}
	\begin{center}
		\subfigure{\includegraphics[width=.25\textwidth]{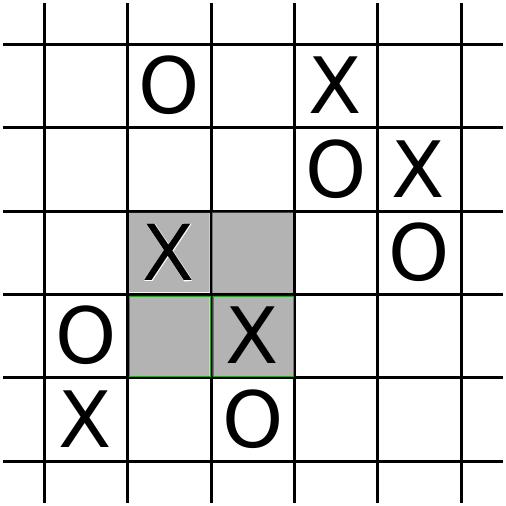}}
		\subfigure{\includegraphics[width=.25\textwidth]{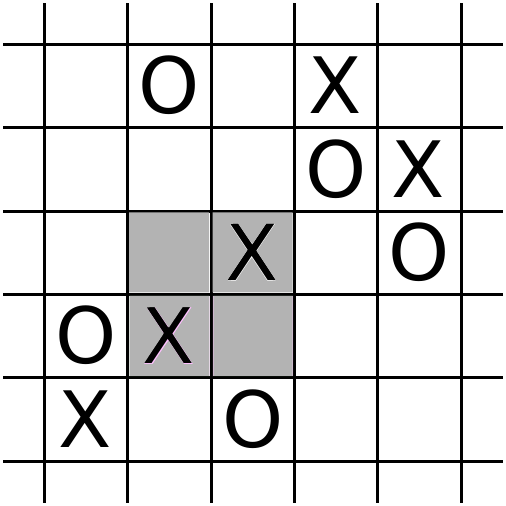}}
	\end{center}
	\caption{A saddle move of the first type on the left-hand grid produces the grid on the right}
	\label{fig:Xsaddle}
\end{figure}

There are two grid moves corresponding to saddles in the cobordism.  The first, an \textbf{$X$-saddle}, which is link-grid move (5) in \cite{MR2915478}, is used when the saddle merges two components of the graph or when it splits one component into two.  If a grid diagram $g$ contains a two-by-two square whose upper left and lower right grid squares contain $X$-markings and whose upper right and lower left grid squares are unoccupied, then doing this saddle move results in a grid diagram $\overline{g}$.  The new grid diagram $\overline{g}$ is exactly the same as $g$ except that in the two-by-two square, the $X$-markings are placed in the upper right and lower left grid squares, with the upper left and lower right squares unoccupied, as shown in \cref{fig:Xsaddle}.

\begin{figure}
	\begin{center}
		\subfigure{\includegraphics[width=.25\textwidth]{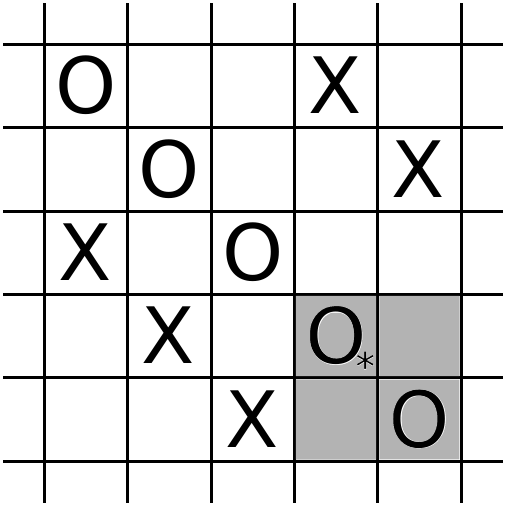}}
		\subfigure{\includegraphics[width=.25\textwidth]{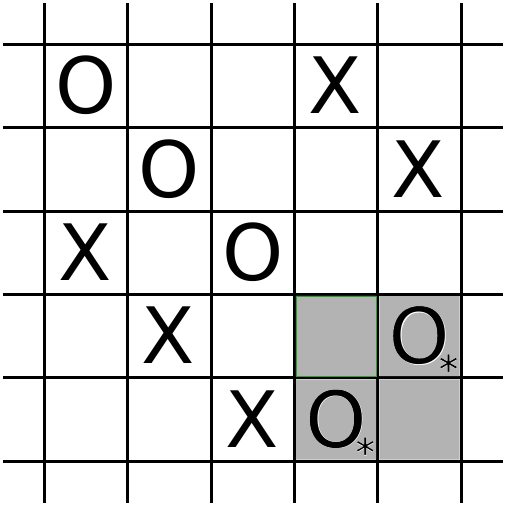}}
	\end{center}
	\caption{A saddle move of the second type on the left-hand grid produces the grid on the right}
	\label{fig:Osaddle}
\end{figure}

The second type of saddle grid move, an \textbf{$O$-saddle}, is used only when the saddle in the cobordism splits one component of the link into two components.  This move is link-grid move (6) from \cite{MR2915478}.  It is exactly the same as the first saddle move except that the two-by-two square which differs in $g$ and $\overline{g}$ contains a special $O$-marking in the upper left corner and a regular $O$ in the lower right corner in $g$, and special $O$-markings in the uper right and lower left corners in $\overline{g}$.  An example is shown in \cref{fig:Osaddle}.

For the proof of the inequality, we will use the combinatorial definition of Alexander grading from \cite{MR2915478}, which we will denote as $A'$.

\begin{definition}
	For a generator $\mathbf{x}$ in a grid diagram $g$, the Alexander grading of $\mathbf{x}$ is 
		\[ A' (\mathbf{x}) = \mathcal{J}(\mathbf{x}, \mathbb{X - O}) - \frac 1 2 \mathcal{J}(\mathbb{X,X}) + \frac 1 2 \mathcal{J} (\mathbb{O,O}) - \frac{n-1}{2},\]
	where $n$ is the grid size of $g$. For an $l$-component link, this definition differs slightly from the usual combinatorial definition of Alexander grading $A(\mathbf{x})$ from \cite{MR2372850}, which can be obtained from $A'(\mathbf{x})$ by adding $\frac{l-1}{2}$.
	\end{definition}

\begin{definition}
	For a tight grid diagram $g$ representing an $l$-component link $L$ in $S^3$, define 		
		\[ \tau' (L) = \text{min} \{ m \in \frac 1 2 \mathbb{Z} | \iota _m \text{ is non-trivial} \}\]
	where  $\widehat{\mathcal{F}}_m'$ is the Alexander filtration induced by the Alexander grading $A'(\cdot)$ and $\iota_m : H_*(\widehat{\mathcal{F}}_m') \rightarrow H_*(\widehat{CF}(g))$ is the map induced by inclusion.
\end{definition}

\begin{lemma}
	The Alexander grading $A(\cdot)$ from \cite{MR2372850} is equal to the Alexander grading defined in \cref{def:Alexgrading}, and for an $l$-component link the $\tau$ defined in \cref{def:tau} is equal to $\tau ' + \frac{l-1}{2}$.
\end{lemma}


\subsection{Link Cobordisms and the $\tau$ Invariant}

\begin{theorem}
If $L_1$ and $L_2$ are $l_1$- and $l_2$-component links, respectively, and $F$ is a connected genus $g$ cobordism from $L_1$ to $L_2$, then 
\[ 1 - g - l_1 \leq \tau (L_1) - \tau (L_2) \leq g + l_2 -1. \]
	\begin{proof}
	The proof will follow the same basic outline of the proof of the main theorem in \cite{MR2915478}.  Consider the cobordism $F$ as a ``movie."  Then there are some number of births, deaths, and saddles in the movie, and the genus $g = \frac 1 2 (s - b - d) +1 - \frac{l_1+l_2}{2}$, where $b,d,$ and $s$ are the number of births, deaths, and saddles, respectively.  We can alter the cobordism slightly so that each of the movie moves happens at a distinct time and so that all of the births take place before any of the saddles, all of the saddles take place before any of the deaths, and the last $l_2 + d - l_1$ saddles split one link component into two.
	
	Note that $l_2 + d - l_1$ is always greater than or equal to $0$. If $l_2 > l_1$, then this is obviously true.  If $l_1 > l_2$, then we must have $d \geq l_1 - l_2$ since both $g_1$ and $g_2$ are tight link grid diagrams and deaths are the only move that reduce the number of special $O$-markings in the grid.  Therefore $d \geq l_1 - l_2 > 0$, so $l_2 + d - l_1 \geq 0$.
	
	As Sarkar shows in \cite{MR2915478}, the modified cobordism can be represented by a sequence of link grid diagrams, such that the first grid, $g_1$ is a tight diagram for $L_1$, the last grid, $g_2$ is a tight diagram for $L_2$, and each diagram in the sequence is obtained from the one before it by a commutation, stabilization, destabilization, birth, $X$-saddle, $O$-saddle, or death grid move, or by renumbering the ordinary $O$-markings.  
	
	As shown in \cite{MR2915478}, the chain maps associated to renumbering the ordinary $O$-markings, commutations, stabilizations, and de-stabilizations are quasi-isomorphisms which preserve both the Maslov and $A'$ Alexander gradings.  The chain maps associated to births is a quasi-isomorphism which preserves the Maslov grading and shift the Alexander grading $A'$ by $- \frac 1 2$.  The chain maps associated to $X$-saddles are the identity maps, and they shift the Alexander grading $A'$ by $+ \frac 1 2$.  The chain maps associated to $O$-saddles induce injective maps on homology and shift the Alexander grading $A'$ by $- \frac 1 2$.  The chain maps associated to deaths induce surjective maps on homology and shift the Alexander grading $A'$ by $+ \frac 1 2$.
	
	Now we will track the overall shift in the Alexander grading $A'$ over the sequence of moves in the (modified) cobordism. Since there are $b$ births, the shift from the births is $- \frac 1 2 b$.  Next we need to figure out how many of the saddles are represented by $X$-saddle grid moves and how many by $O$-saddles.  Any saddle can be represented by either an $X$-saddle or an $O$-saddle grid move, but $O$-saddles and deaths are the only moves that change the number of special $O$-markings in the diagrams. Therefore we can choose which saddles will be represented by $X$-saddles and which by $O$-saddles so that we will have the correct number of special $O$-markings at each stage of the cobordism.  Since the beginning and ending grid diagrams $g_1$ and $g_2$ are tight, we know that $g_2$ has $l_2$ special $O$-markings and $g_1$ has $l_1$ special $O$-markings.  Since the death move removes a special $O$-marking, we need to have $l_2+d$ special $O$-markings after all of the saddles have been performed but before the deaths.  Therefore we should have $l_2+d-l_1$ $O$-saddles in the cobordism, and these are the last saddles.  The fact that we chose that the last $l_2+d-l_1$ saddles should be splits ensures that there will not be more than one special $O$-marking on any one component, so the ending grid diagram $g_2$ will be tight.  The rest of the saddles, which is to say the first $s-l_2-d+l_1$ saddles in the cobordism, are $X$-saddles.
	
	Now we can see that the Alexander grading shift from the $X$-saddles is $+\frac 1 2 (s-l_2-d+l_1)$ and the shift from the $O$-saddles is $-\frac{1}{2} (l_2+d-l_1)$.  Since there are $d$ deaths, the shift from the deaths is $+\frac{1}{2}d$.  Adding up the grading shifts from all of the cobordism moves, the total shift is $\frac{1}{2}(s-b-d) + l_1 - l_2$.
	
	Following \cite{MR2915478}, we know that $\tau '(L_1)$ is less than or equal to $\tau '(L_2)$ plus the Alexander grading shift of the cobordism from $L_1$ to $L_2$.  Therefore 
	\[ \tau '(L_2) \leq \tau ' (L_1) + \frac{1}{2}(s-b-d) + l_1 - l_2,\]
	and after some algebraic manipulation, we see that
	\[ \tau '(L_2) + \frac{l_2-1}{2} - \tau '(L_1) - \frac{l-1}{2} \leq \left( \frac{s-b-d}{2} + 1 + \frac{l_2-l}{2} \right) + l_1 - 1. \]
	Now we observe that since the genus of F is $g = \frac{s-b-d}{2} + 1 + \frac{l_2-l_1}{2}$ and $\tau ' (L_1) + \frac{l_1-1}{2} = \tau (L_1)$, we have
	\[ \tau(L_2) - \tau (L_1) \leq g + l_1 - 1.\]
	
	To prove the other inequality, we reverse the direction of F and consider it as a cobordism from $L_2$ to $L_1$.  Following the same proof as for the first inequality, we see that 
	\[ \tau (L_1) - \tau (L_2) \leq g + l_2 - 1.\]
	
	\end{proof}
	\label{thm:LinkCob}
\end{theorem}


\subsection{Application to link sliceness}

If an $l$-component link $L$ is slice, then there is a concordance between $L$ and the $l$-component unlink.  We can modify this concordance by connect-summing the annuli together and capping off all but one of the unlink's components to produce a connected genus zero cobordism from $L$ to the unknot $U$.  Applying \cref{thm:LinkCob} to this cobordism, we see that 
\[ 1-l \leq \tau(L) - \tau (U) \leq 0. \]
Since $\tau(U) = 0$ we have the following corollary:

\begin{corollary}
\label{cor:linkslice}
	If an $l$-component link $L$ has $\tau (L) > 0$ or $\tau (L) \leq -l$, then $L$ is not slice.
\end{corollary}





\bibliographystyle{amsalpha}
\bibliography{../KatherineBibl.bib}

\end{document}